\documentclass{amsart}

\usepackage[margin=0.75in]{geometry}

\usepackage[utf8]{inputenc}

\usepackage{amssymb, amsmath, amsthm, amsaddr}
\usepackage{array}
\usepackage{bm}
\usepackage{booktabs}
\usepackage{caption}
\usepackage{graphicx}
\usepackage[colorlinks,linkcolor=blue,citecolor=blue]{hyperref}
\usepackage{mdframed}
\usepackage[parfill]{parskip}

\usepackage{xcolor}

\newcommand{\R}{\mathbb{R}}
\newcommand{\T}{\mathcal{T}}
\newcommand{\LOR}{\textit{LOR}}

\makeatletter
\let\div\@undefined
\makeatother
\newcommand{\div}{\nabla \cdot}
\newcommand{\paren}[1]{\left( #1 \right)}
\newcommand{\bracket}[1]{\left\{ #1 \right\}}
\newcommand{\abs}[1]{\left| #1 \right|}
\newcommand{\pd}[2]{\cfrac{\partial #1}{\partial #2}}
\newcommand{\dd}[2]{\cfrac{d #1}{d #2}}
\newcommand{\st}{\; \mid \;}
\renewcommand{\Re}{\mathrm{Re}}

\begin{document}

\title%
[High-order matrix-free incompressible flow solvers]%
{High-order matrix-free incompressible flow solvers with GPU acceleration and low-order refined preconditioners}

\author{Michael Franco}
\address{Department of Mathematics, University of California, Berkeley}
\author{Jean-Sylvain Camier, Julian Andrej, and Will Pazner}
\address{Center for Applied Scientific Computing, Lawrence Livermore National
       Laboratory}

\begin{abstract}
We present a matrix-free flow solver for high-order finite element discretizations of the incompressible Navier-Stokes and Stokes equations with GPU acceleration.
For high polynomial degrees, assembling the matrix for the linear systems resulting from the finite element discretization can be prohibitively expensive, both in terms of computational complexity and memory.
For this reason, it is necessary to develop matrix-free operators and preconditioners, which can be used to efficiently solve these linear systems without access to the matrix entries themselves.
The matrix-free operator evaluations utilize GPU-accelerated sum-factorization techniques to minimize memory movement and maximize throughput.
The preconditioners developed in this work are based on a low-order refined methodology with parallel subspace corrections, as described for diffusion problems in~\cite{Pazner2019Efficient}.
The saddle-point Stokes system is solved using block-preconditioning techniques, which are robust in mesh size, polynomial degree, time step, and viscosity.
For the incompressible Navier-Stokes equations, we make use of projection (fractional step) methods, which require Helmholtz and Poisson solves at each time step.
The performance of our flow solvers is assessed on several benchmark problems in two and three spatial dimensions.
\end{abstract}

\maketitle

\section{Introduction}
\label{sec:introduction}

High-order finite element methods have the potential to attain higher accuracy per degree of freedom than low-order alternatives~\cite{Deville2002,Wang2013}.
Due to these potential computational savings, many high-order methods have been developed to solve a diverse range of computational fluid dynamics (CFD) problems~\cite{Cockburn1989, Bassi1997, Pazner2019High}.
Moreover, the high arithmetic intensity of these algorithms makes them a prime target for use on graphics processing units (GPUs) and GPU-accelerated architectures of current and future supercomputers \cite{Klockner2009,Vermeire2017,Brown2018}.

However, the benefits of high-order methods do not immediately imply that running increasingly higher order simulations for a fixed problem size will result in more efficient simulations. In general, Galerkin finite element methods couple all degrees of freedom (DoFs) within each mesh element.
Therefore, the memory required to store the resulting system matrices grows quadratically with the number of DoFs per element (i.e.~$\mathcal{O}(p^{2d})$ in $d$ spatial dimensions, where $p$ is the polynomial degree).
Furthermore, the naive assembly of the system matrix requires $\mathcal{O}(p^{3d})$ operations, although sum-factorization techniques can lower this to $\mathcal{O}(p^{2d+1})$~\cite{Kronbichler2019}.
Thus, traditional matrix-based approaches are impractical for use with high polynomial degrees, both in terms of computational cost and memory requirements.
Instead, matrix-vector products can be replaced with on-the-fly evaluations of the discretized differential operators in a matrix-free manner.
Combined with sum-factorization techniques on tensor-product elements originally developed in the spectral element community~\cite{Orszag1980, Patera1984}, evaluation of these matrix-free operators can be performed in $\mathcal{O}(p^{d+1})$ operations and $\mathcal{O}(p^d)$ memory~\cite{Pazner2018Approximate}.
Matrix-free evaluations with sum factorization have been shown to outperform sparse matrix-vector products with $p\geq 2$ due to bandwidth bounds on modern architectures~\cite{Kronbichler2012}.
Furthermore, optimized implementations of these operators achieve near peak performance on modern GPUs~\cite{Brown2018,Swirydowicz2019}, making matrix-free high-order operators a desirable choice for performant implementations.

For compressible CFD problems where explicit time stepping can be employed, matrix-free operators implemented on the GPU have been used to great effect~\cite{Klockner2011, Vermeire2017}.
However, most incompressible flow solvers require the solution of large, sparse linear systems~\cite{Deville2002}, thus motivating the development of matrix-free solvers.
Krylov subspace methods are a natural choice for matrix-free solvers, but they require effective preconditioners in order to obtain good performance~\cite{Benzi2005}.
Therefore, in this work we develop matrix-free preconditioners to solve the linear systems arising from high-order tensor-product finite element discretizations of the steady Stokes, unsteady Stokes, and unsteady incompressible Navier-Stokes equations.
Particular emphasis is placed on solver robustness with respect to discretization and mesh parameters.
In recent years, there has been much work on the topic of matrix-free preconditioning for high-order discretizations.
Matrix-free multigrid methods using point Jacobi and Chebyshev smoothing were considered in~\cite{Rudi2015} and~\cite{Kronbichler2019}.
Matrix-free tensor-product approximations to block Jacobi preconditioners for discontinuous Galerkin discretizations were constructed in~\cite{Pazner2018Approximate} and~\cite{Pazner2018b}.
A number of other matrix-free methodologies for high-order discontinuous Galerkin flow solvers have been proposed, using techniques such as multigrid and block preconditioning~\cite{Bastian2019,Franciolini2017,Fehn2018}.
In this work, we extend sparse, low-order refined preconditioners~\cite{Orszag1980,Deville1990,Canuto1985} with parallel subspace corrections, originally described for diffusion problems in~\cite{Pazner2019Efficient}.
The resulting preconditioners are robust in both the mesh size and polynomial degree.

The structure of this paper is as follows.
In \hyperref[sec:equations]{Section~\ref{sec:equations}}, we will briefly describe the equations governing incompressible fluid flow.
\hyperref[sec:spatialDiscretization]{Section~\ref{sec:spatialDiscretization}} will describe the high-order spatial discretization of these equations, with particular emphasis on the implementation of the matrix-free operators. \hyperref[sec:temporalDiscretization]{Section~\ref{sec:temporalDiscretization}} will discuss the temporal discretization, using both Runge-Kutta methods and projection methods.
Matrix-free solvers and preconditioners will be developed in \hyperref[sec:solvers]{Section~\ref{sec:solvers}}.
The performance of our GPU implementation of the matrix-free operators will be analyzed in \hyperref[sec:gpu]{Section~\ref{sec:gpu}}.
Finally, we will present numerical results verifying the $h$- and $p$-robustness of our solvers in 2D and 3D in \hyperref[sec:numericalResults]{Section~\ref{sec:numericalResults}}.

\section{Governing equations of incompressible flow}
\label{sec:equations}

Of interest in this work are the equations governing incompressible fluid flow.
In particular, we consider both the incompressible Stokes and Navier-Stokes equations in $d$ spatial dimensions ($d=2,3$).
The spatial domain is denoted $\Omega \subseteq \R^d$.
The unknowns are velocity $\bm u : \Omega \to \R^d$ and pressure $p : \Omega \to \R$.
The unsteady incompressible Navier-Stokes equations with Dirichlet boundary conditions are
\begin{equation} \begin{aligned}
   \frac{\partial \bm u}{\partial t} + (\bm u \cdot \nabla)\bm u - \nu \Delta \bm u + \nabla p &= \bm{f} &&\quad \text{in }\Omega, \\
   \nabla\cdot\bm u &= 0 &&\quad \text{in }\Omega, \\
   \bm u &= \bm g_D  &&\quad \text{on }\partial\Omega,
\end{aligned} \label{eq:ins} \end{equation}
where we assume a uniform density and viscosity. Here, $\nu$ is the kinematic viscosity, $\bm{f}$ is a known forcing term, and $\bm{g}_D$ is the specified Dirichlet boundary condition.
In the limit of low Reynolds numbers, we are also interested in the linear Stokes equations, which can be obtained by neglecting the nonlinear convection term from~\eqref{eq:ins}, resulting in
\begin{equation}\begin{aligned}
   \frac{\partial \bm u}{\partial t} - \nu \Delta \bm u + \nabla p &= \bm f &&\quad \text{in }\Omega, \\
   \nabla\cdot\bm u &= 0 &&\quad \text{in }\Omega, \\
   \bm{u} &= \bm{g}_D &&\quad \text{on }\partial\Omega .
\end{aligned}\label{eq:td-stokes}\end{equation}
The steady versions of \eqref{eq:ins} and \eqref{eq:td-stokes} will also be considered.

\section{Spatial discretization via high-order finite elements}
\label{sec:spatialDiscretization}

The governing equations are discretized using a high-order continuous finite element method.
First, the spatial domain $\Omega$ is discretized using an unstructured mesh $\T_h$ consisting of tensor-product elements (mapped quadrilaterals in two dimensions and mapped hexahedra in three dimensions).
Note that these element mappings allow for high-order or curved elements.
We define the following finite element function spaces on $\T_h$:
\begin{equation}\begin{aligned}
V_p &= \bracket{\bm{v} \in \paren{H^1(\Omega)}^d \st \bm{v}(K) \in \paren{\mathcal{Q}_p(K)}^d \ \forall K \in \T_h} \\
P_q &= \bracket{s \in H^1(\Omega) \st s(K) \in \mathcal{Q}_q(K) \ \forall K \in \T_h}
\end{aligned}\end{equation}
We will use nodal, Gauss-Lobatto tensor-product bases on each of the spaces $V_p$ and $P_q$.
We write the finite element formulation for \eqref{eq:td-stokes} as: find a velocity-pressure pair $(\bm{u},p) \in (V_p, P_q)$ such that
\begin{equation}\begin{aligned}
\paren{\pd{\bm{u}}{t}, \bm{v}} + \paren{\nu \nabla \bm{u}, \nabla \bm{v}} + \paren{\nabla p, \bm{v}} &= \paren{\bm{f}, \bm{v}} \quad && \forall \bm{v} \in V_p, \\
-\paren{\div \bm{u}, s} &= 0 \quad && \forall s \in P_q.
\label{eq:weak-stokes}\end{aligned}\end{equation}
Expanding out $\bm{u}$, $\bm{v}$, $p$, and $s$ in terms of the bases for their respective spaces, we obtain the semi-discrete Stokes problem:
\begin{equation}\begin{aligned}
\bm{M} \dot{\bm{u}} + \bm{L} \bm{u} + \bm{G} p &= \bm{f}, \\
-D u &= 0,
\end{aligned}\label{eq:unsteadyStokesSemiDiscrete}\end{equation}
where $\bm{u}$ and $p$ are now reused to represent the vectors of coefficients of the high-order polynomial basis functions approximating their continuous counterparts.
$\bm{M}$ is the vector mass matrix, $\bm{L}$ is the vector stiffness matrix, $\bm{G}$ is the gradient operator, and $D$ is the divergence operator with the following definitions:
\begin{align}
\bm{M}_{ij} &= \int_\Omega \bm{\phi}_i \bm{\phi}_j\ dx, \label{eq:vectorMass} \\
\bm{L}_{ij} &= \int_\Omega \nu \nabla \bm{\phi}_i \cdot \nabla \bm{\phi}_j\ dx, \label{eq:vectorDiffusion} \\
\bm{G}_{ij} &= \int_\Omega \bm{\phi}_i \cdot \nabla \psi_j\ dx, \label{eq:gradient} \\
D_{ij} &= \int_\Omega \psi_i \div \bm{\phi}_j\ dx. \label{eq:divergence}
\end{align}
The bases $\bracket{\bm{\phi}_i}_{i=1}^{n_v}$ and $\bracket{\psi_j}_{j=1}^{n_p}$ span the velocity space $V_p$ and pressure space $P_q$, respectively.
In the steady Stokes problem we take $\dot{\bm{u}} = 0$, so \eqref{eq:unsteadyStokesSemiDiscrete} reduces to the linear system
\begin{equation}
\begin{bmatrix} \bm{L} & \bm{G} \\ -D & 0 \end{bmatrix} \begin{bmatrix} \bm{u} \\ p \end{bmatrix} = \begin{bmatrix} \bm{f} \\ 0\end{bmatrix}.
\label{eq:steadyStokesDiscrete}\end{equation}

A similar treatment of the incompressible Navier-Stokes equations \eqref{eq:ins} yields the semi-discrete problem:
\begin{equation}\begin{aligned}
\bm{M} \dot{\bm{u}} + \bm{L} \bm{u} + \bm{N}\paren{\bm{u}} + \bm{G} p &= \bm{f}, \\
-D u &= 0,
\end{aligned}\label{eq:insSemiDiscrete}\end{equation}
where $\bm{N}\paren{\bm{u}}$ is the discretized nonlinear vector-convection term defined by
\begin{equation}
\bm{N}\paren{\bm{u}}_{i} = \int_\Omega \bm{u}^T \paren{\bm{\Phi} \cdot \nabla} \bm{\Phi} \bm{u} \bm{\phi}_i \ dx.
\label{eq:nonlinear}\end{equation}
Equation \eqref{eq:nonlinear} uses the notation $\bm{\Phi}$ to represent the tensor of all $\bm{\phi}_i$, so that $\paren{\bm{\Phi} \cdot \nabla} \bm{\Phi}$ can be viewed as a matrix of size $n_v \times n_v$ where each entry is the vector
\begin{equation}
\left[ \paren{\bm{\Phi} \cdot \nabla} \bm{\Phi}\right]_{ij} = \paren{\bm{\phi}_i \cdot \nabla} \bm{\phi}_j.
\end{equation}

When solving the steady~\eqref{eq:steadyStokesDiscrete} and unsteady Stokes~\eqref{eq:unsteadyStokesSemiDiscrete} problems, we use the Taylor-Hood finite element space, $(\bm{u}, p) \in (V_p, P_{p-1})$, which achieves optimal convergence rates and is stable for orders $p \geq 2$~\cite{Brezzi1991}. When solving the incompressible Navier-Stokes problem~\eqref{eq:ins}, we use the so-called $P_NP_N$ space, $(\bm{u}, p) \in (V_p, P_p)$ \cite{Guermond2006}.

\subsection{Matrix-free operators}
\label{sec:matrixFree}

\begin{figure*}
\centering
\hspace*{\fill}
\includegraphics{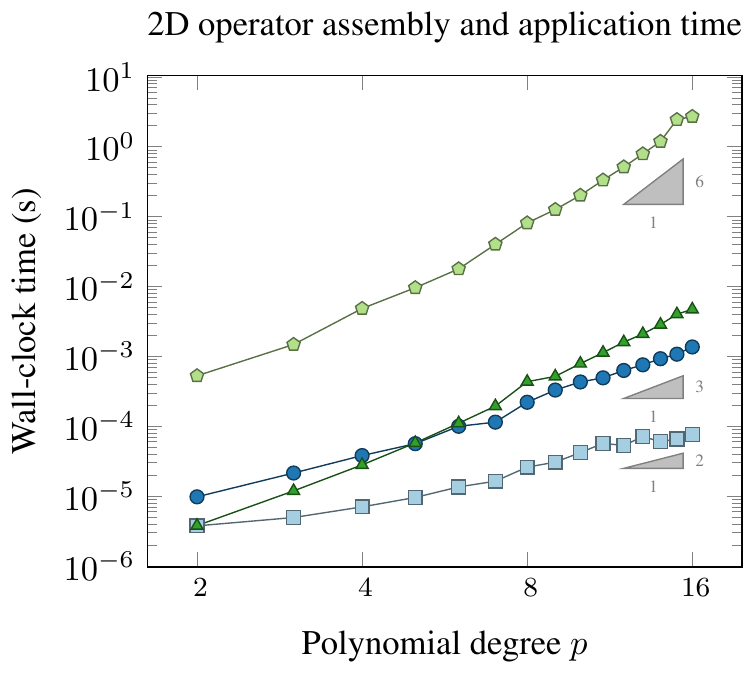}
\hspace*{\fill}
\includegraphics{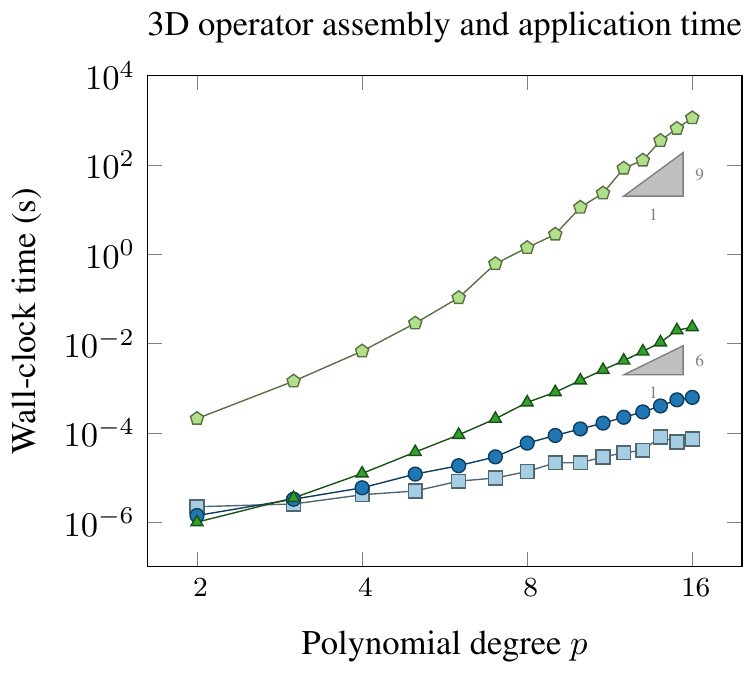}
\hspace*{\fill}

\includegraphics{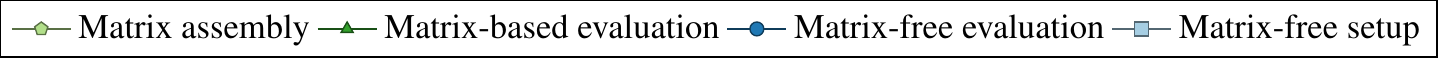}
\caption{Runtime comparison of standard matrix-based and sum-factorized matrix-free operator setup and evaluation for scalar Laplacian in 2D and 3D.
Standard (naive) matrix assembly scales like $\mathcal{O}(p^{3d})$ and matrix-based operator evaluation scales like $\mathcal{O}(p^{2d})$, while matrix-free setup scales like $\mathcal{O}(p^d)$ and matrix-free operator evaluation scales like $\mathcal{O}(p^{d+1})$.}
\label{fig:pa-fa-timings}
\end{figure*}

Explicitly assembling the matrices described by \eqref{eq:vectorMass}--\eqref{eq:divergence} typically costs $\mathcal{O}(p^{3d})$ operations (using the naive algorithm).
Moreover, storing the matrices requires $\mathcal{O}(p^{2d})$ memory.
When running simulations at high polynomial degrees, these costs can be prohibitive.
This is particularly true on GPUs, where the system matrices may be too large to fit in device memory.
Therefore, we turn to sum-factorization techniques to replace our matrix-vector products with on-the-fly evaluations of the integrals~\cite{Swirydowicz2019, Pazner2018Approximate}.
These integrals are efficiently implemented by summing the contributions from each element of the mesh.
For the purposes of illustration, let us consider a matrix-free implementation of the local scalar stiffness matrix on a hexahedral ($d=3$) element $K$:
\begin{equation}
L_{ij}^K = \int_K \nu \nabla \phi_i \cdot \nabla \phi_j\ dx.
\label{eq:localStiffness}\end{equation}
We assume there is an isoparametric mapping $T^K$ from the reference cube $R=[-1,1]^3$ to element $K$ with Jacobian $J^K$. Mapping to the reference cube, we have
\begin{equation}
L_{ij}^K = \int_R \nu \nabla \phi_i^T \paren{J^K}^{-T} \abs{J^K} \paren{J^K}^{-1} \nabla \phi_j\ d\bm{\xi},
\label{eq:localStiffnessMapped}\end{equation}
where now the $\nabla$ operators and basis functions $\bracket{\phi_i}_{i=1}^{n_v}$ are understood to be in the reference space, parameterized by $\bm{\xi} \in \R^d$.
Choosing our bases to be the tensor product of 1D Lagrange interpolating polynomials defined on the $p+1$ Gauss-Lobatto points, we have
\begin{equation}
\phi_{ijk}(\bm{\xi}) = \phi_i(\xi_1)\phi_j(\xi_2)\phi_k(\xi_3), \end{equation}
for all $0 \leq i,j,k \leq p$.
Choosing $n_q^d$ tensor product quadrature points to evaluate the integrals in \eqref{eq:localStiffnessMapped}, we can also denote the quadrature nodes and weights using this same multi-index notation.
That is, the quadrature weights and points are $\bracket{w_{i_qj_qk_q}}_{i_q,j_q,k_q=1}^{n_q}$ and $\bracket{\bm{\xi}_{i_qj_qk_q}}_{i_q,j_q,k_q=1}^{n_q}$.
Thus, \eqref{eq:localStiffnessMapped} becomes
\begin{multline}
L_{ij}^K = \sum_{i_q,j_q,k_q=1}^{n_q} \nu w_{i_qj_qk_q} \nabla \phi_i^T\paren{\bm{\xi}_{i_qj_qk_q}} \paren{J^K(\bm{\xi}_{i_qj_qk_q})}^{-T} \\ \abs{J^K(\bm{\xi}_{i_qj_qk_q})}
\paren{J^K(\bm{\xi}_{i_qj_qk_q})}^{{-1}} \nabla \phi_j\paren{\bm{\xi}_{i_qj_qk_q}}.
\label{eq:localStiffnessQuadrature}\end{multline}

As an aside, the solution $u$ evaluated at a quadrature point is
\begin{align}
u\paren{\bm{\xi}_{i_qj_qk_q}} &= \sum_{i,j,k=0}^{p} u_{ijk} \phi_{ijk}\paren{\bm{\xi}_{i_qj_qk_q}} \\
&= \sum_{k=0}^{p} \phi_k(\xi_{k_q}) \sum_{j=0}^{p} \phi_j(\xi_{j_q}) \sum_{i=0}^{p} u_{ijk} \phi_i(\xi_{i_q}). \label{eq:sumFactorized}
\end{align}
We immediately notice that the number of operations required to evaluate a function at each of the quadrature points is $\mathcal{O}(p^4)$, or $\mathcal{O}(p^{d+1})$ in general.

Kronecker products allow us to simplify notation.
First, we define the one-dimensional Gauss point evaluation matrix as the $n_q \times (p+1)$ Vandermonde-type matrix obtained by evaluating each of the one-dimensional basis functions at all of the quadrature points:
\begin{equation}
B_{i_q,j} = \phi_j(\xi_{i_q}).
\label{eq:B1D}\end{equation}
With this notation, $\eqref{eq:sumFactorized}$ is equivalent to the computation of the Kronecker product
\begin{equation}
u\paren{\bm{\xi}} = \paren{B \otimes B \otimes B} u.
\end{equation}
Likewise, we define the 1D Gauss point differentiation matrix as
\begin{equation}
D_{i_q,j} = \dd{\phi_j}{\xi} \paren{\xi_{i_q}}.
\label{eq:D1D}\end{equation}

Returning to \eqref{eq:localStiffnessQuadrature}, we are now ready to recast the operator as a Kronecker product of local 1D matrices.
First, we precompute the tensor $W \in \R^{n_q^d \times n_q^d \times d \times d}$, which is diagonal in its first two dimensions and is defined by
\begin{equation}
W_{i_qj_qk_q,i_qj_qk_q} = \nu w_{i_qj_qk_q} \paren{J^K(\bm{\xi}_{i_qj_qk_q})}^{-T} \abs{J^K(\bm{\xi}_{i_qj_qk_q})} \paren{J^K(\bm{\xi}_{i_qj_qk_q})}^{-1}.
\end{equation}
Next, we recognize that evaluating $\nabla \Phi$ at all the quadrature points means contracting with the tensor $\bm{G}_\phi \in \R^{n_q^d \times (p+1)^d \times d}$ defined by
\begin{equation}
\bm{G}_\phi = \begin{bmatrix}B \otimes B \otimes D \\
B \otimes D \otimes B \\
D \otimes B \otimes B \end{bmatrix}.
\label{eq:tensorGradient}\end{equation}
Therefore, the local operator \eqref{eq:localStiffnessQuadrature} is equivalent to
\begin{equation}
L^K = \bm{G}_\phi^T W \bm{G}_\phi,
\end{equation}
where the transpose is taken only over the first two dimensions of $\bm{G}_\phi$.
A major benefit of reformulating the local operator into this form is that it is no longer necessary to form a global matrix.
Instead, the action of this operator can be recreated using the 1D matrices $B$ and $D$ and the precomputed tensor $W$.
Precomputing $W$ requires $\mathcal{O}(p^d)$ operations and $\mathcal{O}(p^d)$ memory.
Applying the operator $\bm{G}_\phi$ and its transpose requires $\mathcal{O}(p^{d+1})$ operations and $\mathcal{O}(p^d)$ memory.
If $W$ is computed by interpolating from the nodal points using sum-factorization techniques, then the setup costs would increase to $\mathcal{O}(p^{d+1})$.
To achieve our leading-order costs, we assume that $J$ (and therefore $W$) can be computed or is available at quadrature points in $\mathcal{O}(1)$ time.
Indeed, our implementation stores $J$ at quadrature points, as described in the CEED framework \cite{Dobrev2017}.

Using these sum-factorization techniques, we can also derive matrix-free implementations of $\bm{M}$, $\bm{G}$, $D$, and $\bm{N}(\bm{u})$ from their definitions.
Each operator has comparable computational cost and memory requirements as in the case of this scalar diffusion operator.

Figure \ref{fig:pa-fa-timings} shows run times for sum-factorized matrix-free operator setup and evaluation compared with standard matrix assembly and matrix-based operator evaluation for polynomial degrees between 2 and 16.
We notice that matrix assembly is the most expensive operation for all the cases tested, usually taking about two orders of magnitude more time than operator evaluation.
In 2D, our implementation of matrix-free operator evaluation is more efficient than matrix-based operator evaluation for polynomial degrees greater than 5.
In 3D, the matrix-free sum-factorized operator evaluation is more efficient for polynomial degrees greater than 2.
In the matrix-free context, the setup and precomputations typically represent a negligible portion of the overall cost of using these operators.

\section{Temporal discretization}
\label{sec:temporalDiscretization}

For the time dependent problems~\eqref{eq:ins} and~\eqref{eq:td-stokes}, we discretize in time using several time integration methods.
Broadly speaking, these are classified as split (i.e.~projection or fractional step) methods, or unsplit methods.
For the unsplit methods, we use the method of lines to first discretize in space and then temporally discretize the resulting system of ordinary differential equations \eqref{eq:insSemiDiscrete} and \eqref{eq:unsteadyStokesSemiDiscrete}.
Here, we use diagonally implicit Runge-Kutta (DIRK) schemes as our time-integration method~\cite{Alexander1977}.
On the other hand, split methods such as projection-type methods can be developed in order to decouple the solution of the velocity and pressure components.
As a result, these methods can be computationally efficient, at the cost of incurring splitting and other approximation errors.
Each of these methods will be discussed in greater detail in the following sections.

\subsection{DIRK methods}
\label{sec:dirk}

Consider the system of ordinary differential equations
\begin{equation}
   M \dot{\bm{y}} = \bm r(t, \bm y),
\end{equation}
obtained from a finite element spatial discretization.
Here, $\bm y$ represents a vector of degrees of freedom, $M$ the finite element mass matrix, and $\bm r$ the potentially nonlinear finite element residual vector.
Let $\bm y^n$ denote the known solution at time $t^n$.
A general $s$-stage Runge-Kutta method to approximate the solution at time $t^{n+1} = t^n + \Delta t$ can be written as
\begin{align}
   \label{eq:rk-k}
   M \bm{k}_i^n &= \bm{r}\left( t^n + \Delta t c_i,
                  \bm{y}^n + \Delta t \sum_{j=1}^s a_{ij} \bm{k}_j^n \right), \\
   \label{eq:rk-y1}
   \bm{y}^{n+1} &= \bm{y}^n + \Delta t \sum_{i=1}^s b_i \bm{k}_i^n,
\end{align}
where the coefficients $a_{ij}, b_i,$ and $c_i$ can be expressed compactly in
the form of the \textit{Butcher tableau},
\begin{equation}
   \begin{array}{c|ccc}
       c_1    & a_{11} & \cdots & a_{1s} \\
       \vdots & \vdots & \ddots & \vdots \\
       c_s    & a_{s1} & \cdots & a_{ss} \\
       \hline
              & b_1    & \cdots & b_s
   \end{array}
   =
   {\setlength\extrarowheight{5pt}
   \begin{array}{c|c}
       \bm{c} & A\\
       \hline
       & \bm{b}^T
   \end{array}}.
\end{equation}
A Runge-Kutta scheme is diagonally implicit if its Butcher matrix $A$ is lower triangular, allowing for the solution of the system of equations \eqref{eq:rk-k} through a forward-substitution procedure.
Of particular interest are the so-called \textit{singly} diagonally implicit Runge-Kutta schemes, which use identical diagonal coefficients $a_{ii} = \alpha$.
These schemes allow for reuse of preconditioners, since each stage requires the solution of a backward Euler system with the same time step $\alpha\Delta t$.

The application of DIRK methods to the systems \eqref{eq:ins} and \eqref{eq:td-stokes} is complicated by the fact that there is no temporal evolution equation corresponding to the pressure.
Thus, after spatial discretization we obtain \eqref{eq:insSemiDiscrete} and \eqref{eq:unsteadyStokesSemiDiscrete}, which are systems of diff\-erential-alge\-braic equations (DAEs) \cite{John2006,Rang2007}.
To avoid this difficulty, we reformulate the DIRK method in a manner suitable for solving DAEs \cite{Hairer1996,Brenan1995}.
We first require that the DIRK scheme be \textit{stiffly accurate}, i.e.~that $a_{si} = b_i$ and $c_s = 1$.
Then, consider the general differential-algebraic system
\begin{align}
      M \dot{\bm{y}} &= \bm r(t, \bm y, \bm z), \\
      0 &= \bm g(t, \bm y).
\end{align}
We define the approximate solution to the differential variable $\bm y$ at the $i$th stage by
\begin{equation} \label{eq:rk-y}
   \bm y_i^n = \bm y^n + \Delta t \sum_{j=1}^i a_{ij} \bm{k}_j^n.
\end{equation}
Analogously to \eqref{eq:rk-k}, the stage derivatives $\bm{k}_i^n$ are given by
\begin{align}
   \label{eq:rk-k-dae}
   M \bm{k}_i^n &= \bm{r}( t^n + \Delta t c_i, \bm y_i^n, \bm z_i^n),
\end{align}
Multiplying \eqref{eq:rk-y} by the mass matrix $M$ and inserting \eqref{eq:rk-k-dae}, we obtain
\begin{equation} \label{eq:rk-dae}
   M \bm y_i^n = M \bm y^n + \Delta t \sum_{j=1}^i a_{ij} \bm{r}( t^n + \Delta t c_j, \bm y_j^n, \bm z_j^n),
\end{equation}
which, when augmented with the constraint
\begin{equation}
   0 = \bm g(t, \bm y_i),
\end{equation}
results in a system of equations for the $i$th stage approximations for the differential and algebraic variables $\bm y_i$ and $\bm z_i$.
Because the DIRK schemes under consideration are stiffly accurate, the values at the next time step are given by the final stage approximations
\begin{equation}
   \bm y^{n+1} = \bm y_s^n, \qquad
   \bm z^{n+1} = \bm z_s^n.
\end{equation}

Applying this DIRK method to the semi-discrete Stokes problem~\eqref{eq:unsteadyStokesSemiDiscrete} requires solving the linear system
\begin{equation}
\begin{bmatrix} \frac{1}{\alpha \Delta t}\bm{M}+\bm{L} & \bm{G} \\ -D & 0 \end{bmatrix} \begin{bmatrix} \bm{u}_i^n \\ p_i^n \end{bmatrix} = \begin{bmatrix} \bm{F}_i \\ 0\end{bmatrix}
\label{eq:unsteadyStokesDiscrete}\end{equation}
every Runge-Kutta stage, with right-hand side $\bm{F}_i$ given by
\begin{equation}
\bm{F}_i = \cfrac{\bm{M}\bm{u}^n}{\alpha \Delta t} + \bm{f}_i + \cfrac{1}{\alpha} \sum_{j=1}^{i-1} a_{ij}(\bm{f}_j-\bm{L}\bm{u}_j^n - Gp_j^n).
\end{equation}
This fully discrete linear system~\eqref{eq:unsteadyStokesDiscrete} and its steady state counterpart~\eqref{eq:steadyStokesDiscrete} are saddle point problems~\cite{Benzi2005}.
We present robust solvers for such saddle point systems in \hyperref[sec:blockPrecond]{Section~\ref{sec:blockPrecond}}.

\subsection{Projection methods}
\label{sec:projection}

Projection methods, first introduced by Chorin in 1967, are a class of split  methods for the temporal integration of the incompressible Navier-Stokes equations~\cite{Chorin1967,Chorin1968}.
These methods have the attractive feature that they only require the solution to uncoupled, positive-definite problems, instead of the coupled, saddle-point type problems that are required by the DIRK methods described in \hyperref[sec:dirk]{Section~\ref{sec:dirk}}.
For this reason, projection and fractional-step methods have become extensively used for incompressible flow problems~\cite{Deville2002,Prohl1997}.
The original method of Chorin has since been modified and extended to a wide range of variants~\cite{Bell1989,Brown2001,Karniadakis1991,Tomboulides1997,Orszag1986}.
See e.g.~\cite{Guermond2006} for a review and analysis of a selection of these variants.

Following the method presented in~\cite{Tomboulides1997}, we use equal order polynomial degrees for velocity and pressure, often known as a $P_N P_N$ formulation.
This method uses an implicit-explicit time-integration scheme for the viscous and convective terms respectively, thereby avoiding the need to solve a nonlinear system of equations at every time step.
We use a BDF method of order $k$ for the implicit terms and an extrapolation method of order $k$ for the explicit terms with corresponding coefficients $b_j$ and $a_j$~\cite{Tomboulides1997,Orszag1986}.
First, we introduce the linear term $L(\bm u) = \nu \Delta \bm u$ and nonlinear term $N(\bm u) = -(\bm u \cdot \nabla) \bm u$ as well as their time-extrapolated versions,
\begin{align}
    L^*(\bm u^{n+1}) &= \sum_{j=1}^{k} a_j L(\bm u^{n+1-j}),\\
    N^*(\bm u^{n+1}) &= \sum_{j=1}^{k} a_j N(\bm u^{n+1-j}).
\end{align}
Directly applying a BDF method to~\eqref{eq:ins} yields
\begin{equation}
\sum_{j=0}^{k} \frac{b_j}{\Delta t} \bm u^{n+1-j} = -\nabla p^{n+1} + L(\bm u^{n+1}) + N^*(\bm u^{n+1}) + \bm f^{n+1},
\label{eq:projtdiscfull}\end{equation}
where $\bm{f}^{n+1}$ is assumed to be known \textit{a priori}.
Introducing $\bm{F}^*(\bm{u}^n)$ to represent all known terms at a given time step,
\begin{equation}
\bm{F}^*(\bm{u}^n) = - \sum_{j=1}^{k} \frac{b_j}{\Delta t} \bm u^{n+1-j} + N^*(\bm u^{n+1}) + \bm f^{n+1},
\label{eq:projKnown}\end{equation}
we can simplify \eqref{eq:projtdiscfull} to
\begin{equation}
\frac{b_0}{\Delta t} \bm u^{n+1} = -\nabla p^{n+1} + L(\bm u^{n+1}) + \bm F^*(\bm u^n).
\label{eq:projtdiscexpl}\end{equation}
Unfortunately, despite using a $k > 1$ order time-integration scheme, this method yields at most first-order convergence in time for velocity, as shown in~\cite{Karniadakis1991} and later proved by~\cite{Guermond2006}.
This is caused by splitting errors and large divergence errors on the boundary of the domain.
Therefore we use the velocity-correction formulation presented in~\cite{Karniadakis1991}, where the linear term $L(\bm u)$ is instead expressed as
\begin{equation}
L_{\times}(\bm u) = \nu \nabla (\nabla \cdot \bm u) - \nu \nabla \times \nabla \times \bm u,
\label{eq:projLinear}\end{equation}
using well-known vector calculus identities.
This alternative form of the linear term imposes the incompressibility constraint from~\eqref{eq:ins} weakly, by setting the first term in~\eqref{eq:projLinear} equal to zero.

In order to solve for pressure, we first rearrange \eqref{eq:projtdiscexpl} and take the divergence of both sides in order to get
\begin{align}
\nabla p^{n+1} &= -\frac{b_0}{\Delta t} \bm{u}^{n+1} + \underbrace{L_{\times}^*(\bm u^{n+1}) + F^*(\bm u^{n+1})}_{\tilde{F}^*(\bm u^{n+1})}, \label{eq:projpresreord} \\
\implies \Delta p^{n+1} &= \nabla \cdot \tilde{F}^*(\bm u^{n+1}). \label{eq:pres-poisson}
\end{align}
Notice how the first right hand side term of \eqref{eq:projpresreord} vanishes due to the incompressibility constraint.
Equation~\ref{eq:pres-poisson} is closed by the boundary condition,
\begin{equation}
\nabla p^{n+1} \cdot \hat{\bm{n}} = - \frac{b_0}{\Delta t} \bm u^{n+1} \cdot \hat{\bm{n}} + \tilde{F}^*(\bm u^{n+1}) \cdot \hat{\bm{n}} \quad \text{on }\partial\Omega,
\label{eq:presPoissonBC}\end{equation}
where $\hat{\bm{n}}$ is the outward pointing normal vector.
We use the known Dirichlet boundary condition $\bm{u}^{n+1} \cdot \hat{\bm{n}} = \bm{g}_D^{n+1} \cdot \hat{\bm{n}}$ to evaluate~\eqref{eq:presPoissonBC}.
In the case of a pure Neumann boundary condition, we close the system with a mean-zero condition on pressure:
\begin{equation}
\int_{\Omega} p\ dx = 0.
\label{eq:meanZero}\end{equation}

Therefore, this projection method computes $\bm u^{n+1}$ in three steps.
First, the extrapolated contributions from the nonlinear and forcing terms are combined to compute $\bm  F^*(\bm u^n)$ via~\eqref{eq:projKnown}.
Second, we solve for $p^{n+1}$ in the pressure-Poisson problem \eqref{eq:pres-poisson}, closed with \eqref{eq:presPoissonBC} and \eqref{eq:meanZero}.
Finally, we return to~\eqref{eq:projtdiscexpl} and solve for $\bm u^{n+1}$ in the following Helmholtz problem:
\begin{align} \label{eq:proj-helmholtz}
\frac{b_0}{\Delta t} \bm u^{n+1} - L(\bm u^{n+1}) &= -\nabla p^{n+1} + \bm F^*(\bm u^{n+1}) && \text{in }\Omega,\\
\bm u^{n+1} &= \bm g_D^{n+1} && \text{on }\partial\Omega.
\end{align}

This projection method is $k$th order in time for velocity (up to $k=3$)~\cite{Guermond2006}.
As previously mentioned, a major benefit of this method is its computational efficiency.
Each time step requires only one new nonlinear evaluation, one Poisson solve, and one Helmholtz solve.
We will discuss the matrix-free solvers that we use for these sub-problems in \hyperref[sec:projectionPrecond]{Section~\ref{sec:projectionPrecond}}.

\section{Solvers and matrix-free preconditioners}
\label{sec:solvers}

The numerical methods described above require solving large, sparse linear systems.
The fully discrete, steady Stokes equation requires solving the saddle-point linear system~\eqref{eq:steadyStokesDiscrete}.
The time discretization of the unsteady Stokes equation \eqref{eq:td-stokes} by a DIRK method results in a sequence of saddle point problems~\eqref{eq:unsteadyStokesDiscrete}.
The velocity-correction schemes require the solution of a Poisson problem~\eqref{eq:pres-poisson} for the pressure and a Helmholtz equation~\eqref{eq:proj-helmholtz} for the velocity.
Additionally, the nonlinear extrapolation requires the inversion of the velocity mass matrix.

In order to solve these systems, we make use of preconditioned Krylov subspace methods, such as the conjugate gradient (CG) and generalized minimal residual (GMRES) methods~\cite{Saad2003}.
These iterative solvers are a natural choice for matrix-free solvers, since they only require the action of the operator, which we compute using the matrix-free algorithms described in \hyperref[sec:matrixFree]{Section~\ref{sec:matrixFree}}, and the evaluation of a preconditioner.
The main challenge associated with the matrix-free solution of high-order flow problems is constructing efficient preconditioners that result in iteration counts that are independent of the discretization parameters $h$, $p$, and $\Delta t$.
In this section, we describe the construction of robust preconditioners that do not require the assembly of the high-order system matrices.
We begin by describing our preconditioning strategy for the relatively simpler sub-problems, which can then be combined to create effective preconditioners for the more challenging, coupled problems.

\subsection{Collocated mass preconditioning}
\label{sec:massPrecond}

In order to precondition the mass matrix, we make use of a diagonal preconditioner based on collocated quadrature~\cite{Ewing1990}.
Because we use a nodal Gauss-Lobatto basis for the finite element spaces, integrating at the same set of quadrature points results in a diagonal matrix, which is constructed and inverted in constant time and memory per degree of freedom.
This matrix is spectrally equivalent to the fully-integrated mass matrix, with constants of equivalence independent of $h$ and $p$, and so the number of solver iterations remains uniformly bounded with respect to the mesh size and polynomial degree~\cite{Teukolsky2015}.

\subsection{Low-order refined preconditioners for Poisson and Helm\-holtz problems}
\label{sec:poissonPrecond}

Matrix-free preconditioners for the symmetric positive definite Poisson and Helmholtz problems form the fundamental building blocks for our robust fluid solvers.
These preconditioners are described in detail in \cite{Pazner2019Efficient}, and are based on the spectral equivalence between the high-order finite element discretization, and a low-order ($p_{\textit{low}}=1$) finite element discretization on a Gauss-Lobatto refined mesh.
This equivalence is often refereed to as the finite element method--spectral element method (FEM-SEM) equivalence \cite{Canuto2007}.
The low-order finite element discretization results in a sparse matrix with $\mathcal{O}(1)$ nonzeros per row independent of $p$, the polynomial degree of the original high-order discretization.
Therefore, the memory requirements and computational cost to assemble the low-order matrix are both optimal, scaling linearly in the number of degrees of freedom.

\begin{figure}
\centering
\includegraphics[width=0.38\linewidth]{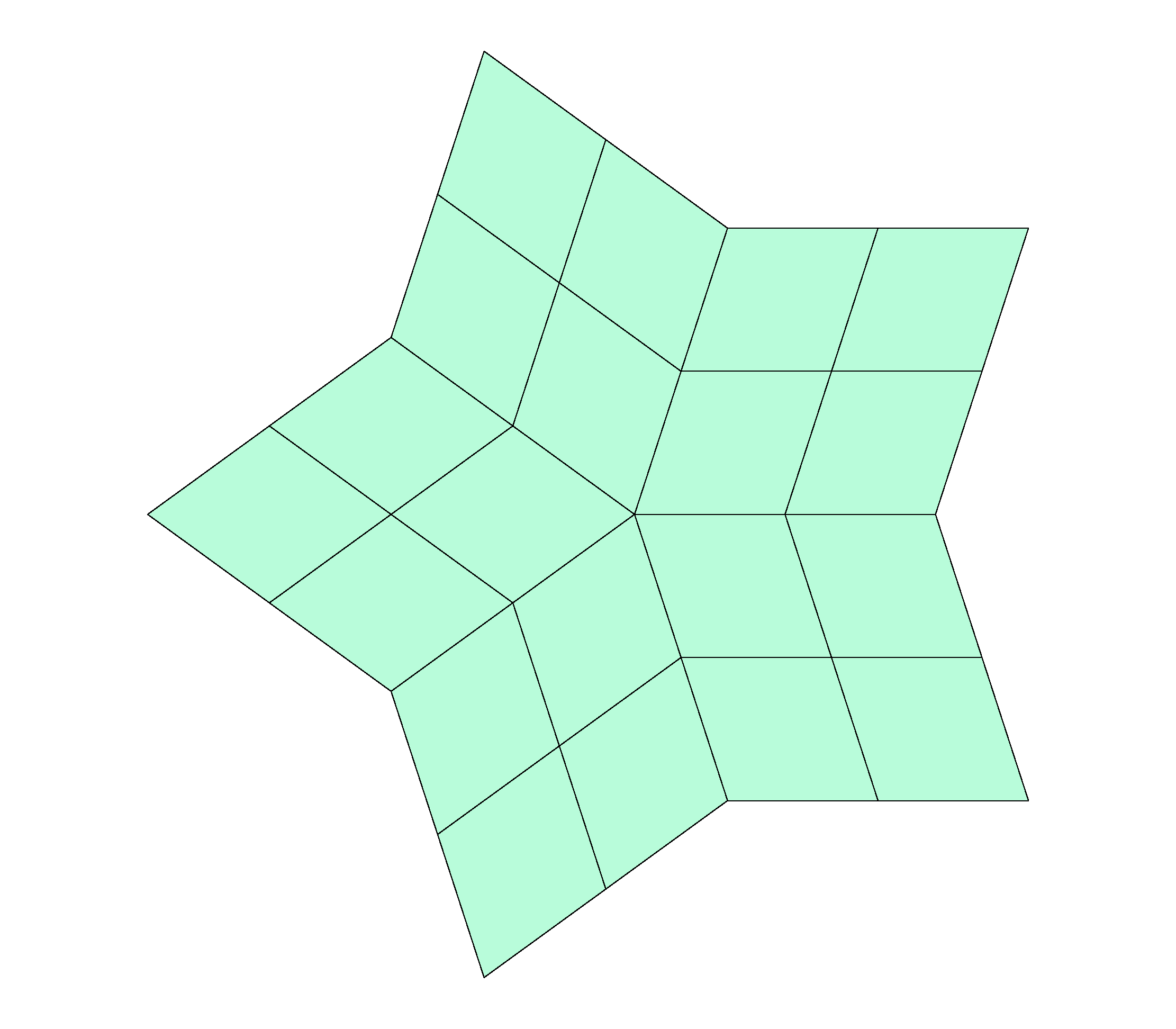}
\includegraphics[width=0.38\linewidth]{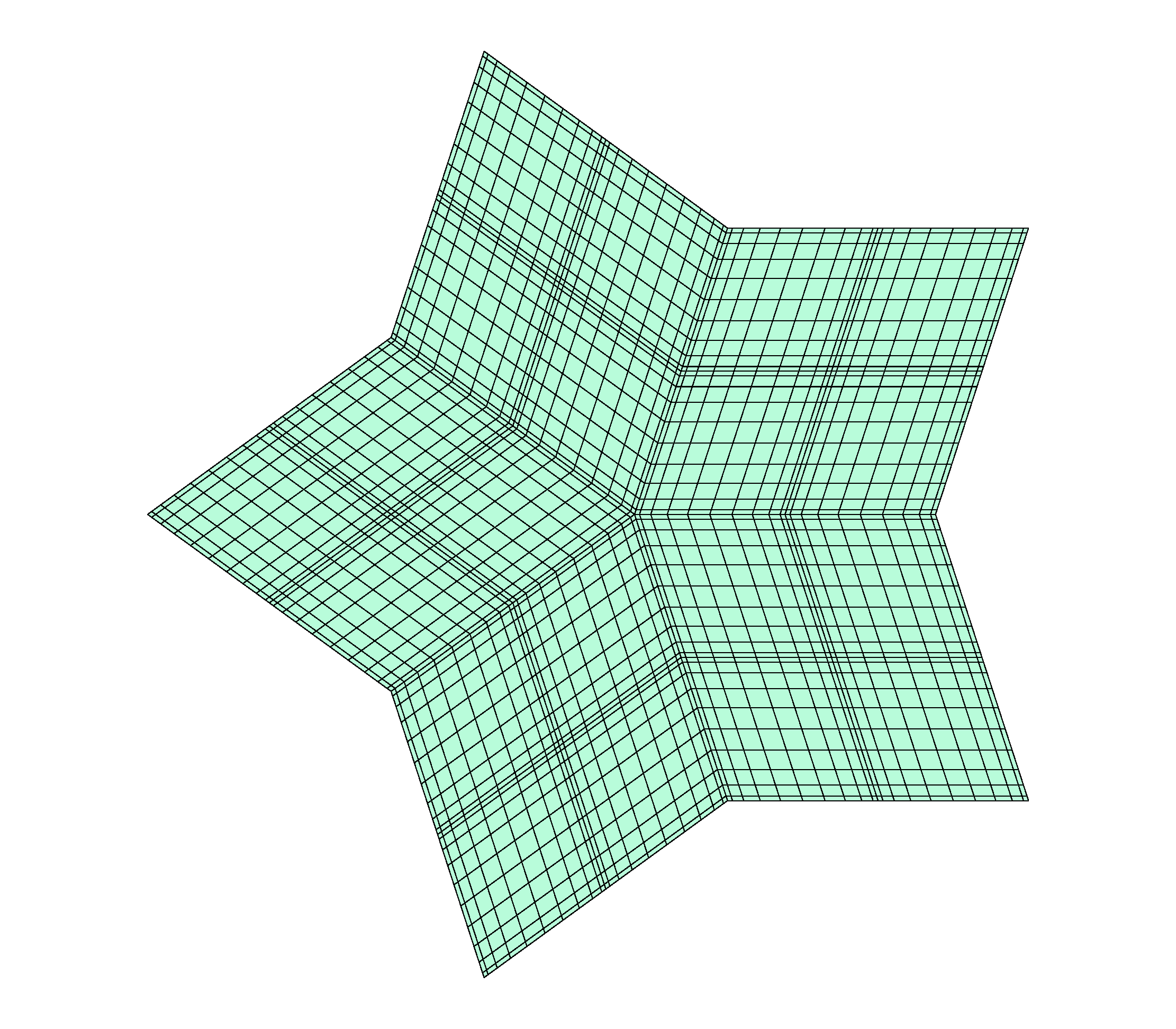}
\caption{Illustration of the low-order refined methodology with $p=10$, showing high aspect ratio elements near the coarse element interfaces.
Left: original high-order mesh $\T_h$. Right: Gauss-Lobatto refined mesh $\T_{\LOR}$.}
\label{fig:lor-mesh}
\end{figure}

Consider a scalar Poisson or Helmholtz problem
\begin{equation}
   A \bm u = \bm b,
\end{equation}
for $A = c M + L$, where $c$ is a non-negative (but possibly zero) constant.
We begin by constructing a \textit{low-order refined} (LOR) operator $A_{\LOR} = c M_{\LOR} + L_{\LOR}$.
Each of the LOR operators is obtained by a standard low-order finite element discretization on a refined mesh $\T_{\LOR}$.
This mesh is obtained by subdividing each element $K \in \T_h$ into the parallelepipeds defined by the Cartesian product of the $p+1$ one-dimensional Gauss-Lobatto points.
Figure~\ref{fig:lor-mesh} illustrates one such low-order refined mesh.
The identity operator can be used to map DoFs from the high-order finite element space using a Gauss-Lobatto basis to the low-order refined space.
It can be shown that the low-order refined mass matrices and stiffness matrices $M_{\LOR}$ and $L_{\LOR}$ are spectrally equivalent to their high-order counterparts, $M$ and $L$ \cite{Canuto1994,Canuto2007,Canuto2010}.
Therefore, $A_{\LOR}$ is spectrally-equivalent to $A$, and a robust preconditioner for $A_{\LOR}$ is, in turn, a robust preconditioner for the original high-order matrix $A$.
The advantage of the matrix $A_{\LOR}$ over $A$ is its greatly increased sparsity, requiring only $\mathcal{O}(1)$ nonzeros per row.
As a consequence, this matrix can be explicitly assembled and stored, allowing for the construction of sophisticated preconditioners.
Note that it is also possible to decompose each coarse element $K \in \T_h$ into a simplicial submesh using the Gauss-Lobatto points.
It has been shown that preconditioners resulting from finite element discretizations using simplicial decompositions can result in improved convergence when compared with the tensor-product decomposition used in this work \cite{Fischer1997,Canuto2010}.
Several new configurations for the low-order mesh were also considered in \cite{Bello-Maldonado2019}.
However, the tensor-product decomposition has the advantage that the implementation may be simplified by reusing the discretization primitives in both the high-order and low-order methods.

The main challenge associated with constructing effective preconditioners for $A_{\LOR}$ is the high aspect ratio associated with the low-order refined mesh $\T_{\LOR}$~\cite{Lottes2005}.
Because the Gauss-Lobatto points are clustered near the endpoints of the interval, the resulting Cartesian product mesh consists of parallelepipeds with aspect ratios that scale like $p$ \cite{Brix2014}.
As a result, the mesh $\T_{\LOR}$ is not shape-regular with respect to $p$, and standard multigrid-type methods will not result in uniform convergence under $p$-refinement.
In order to address this issue, we make use of a structured geometric multigrid V-cycle with ordered ILU smoothing to treat the anisotropy of the problem.
This ordering ensures the number of multigrid iterations is independent of $p$.
In order to obtain an algorithm that is scalable in parallel, we use an overlapping patch-based additive Schwarz method, that is shown to be robust in $h$ and $p$~\cite{Pazner2019Efficient}.
The coarse mesh $\T_h$ is decomposed into overlapping, unstructured vertex patches, together with a low-order coarse space.
The multigrid algorithm with ILU smoothing described above is applied in parallel to each of the subdomains independently.
The additive Schwarz preconditioner is obtained by taking the sum of the subdomain contributions.
A detailed description and analysis of this algorithm is presented in~\cite{Pazner2019Efficient}.

\subsection{Application to the projection method}
\label{sec:projectionPrecond}

The projection method described in \hyperref[sec:projection]{Section~\ref{sec:projection}} requires solvers for the vector mass matrix, the pressure Poisson problem~\eqref{eq:pres-poisson} and the Helmholtz problem~\eqref{eq:proj-helmholtz}.
Each of these problems is symmetric positive-definite, and so we can use a preconditioned conjugate gradient solver with the preconditioners described in \hyperref[sec:massPrecond]{Section~\ref{sec:massPrecond}} and \hyperref[sec:poissonPrecond]{Section~\ref{sec:poissonPrecond}}.
Because the pressure Poisson problem~\eqref{eq:pres-poisson} with pure Neumann conditions has a null space consisting of constant functions, care must be taken to ensure that the right-hand side is orthogonal to the null space.
Therefore, each application of the preconditioner is augmented with an orthogonalization step to ensure convergence.

\subsection{Block preconditioners for Stokes}
\label{sec:blockPrecond}

In order to solve the coupled Stokes problems \eqref{eq:steadyStokesDiscrete} and \eqref{eq:unsteadyStokesDiscrete}, we make use of block-preconditioning techniques, allowing us to reuse the preconditioners for the Poisson and Helmholtz problems.
Consider the representative saddle point system,
\begin{equation}
\begin{bmatrix} \bm{A} & \bm G \\ D & 0 \end{bmatrix} \begin{bmatrix} \bm{u} \\ p \end{bmatrix} = \begin{bmatrix} \bm{f} \\ 0\end{bmatrix}.
\label{eq:saddlePoint}\end{equation}
We define the block triangular preconditioner
\begin{equation}
P_t = \begin{bmatrix} \bm{A} & 0 \\ D & S\end{bmatrix},
\label{eq:blockTriangularPreconditioner}
\end{equation}
where $S = -D\bm A^{-1}\bm G$ is the Schur complement.
Because the preconditioner $P_t$ is not symmetric, we use the flexible GMRES method~\cite{Saad1993}, instead of CG, to solve the preconditioned system.
It can be shown that the spectrum of the preconditioned system $P_t^{-1}\bm A$ consists only of a single eigenvalue, and hence GMRES will converge in at most two iterations~\cite{Benzi2005}.
Applying this preconditioner requires solving linear systems with the matrices $\bm A$ and $S$.
However, the Schur complement $S$ is in general dense, and hence is impractical to form or invert.
Therefore, in this work, we provide a matrix-free operator $\tilde{S}^{-1}$ that approximates the action of $S^{-1}$.

Additionally, we replace the action of $\bm A^{-1}$ with the application of several uniformly-preconditioned conjugate gradient iterations on $\bm A$.
For the steady Stokes system, $\bm{A} = \bm{L}$, so we can use the matrix-free low-order refined preconditioner defined in \hyperref[sec:poissonPrecond]{Section~\ref{sec:poissonPrecond}}.
Likewise, for the unsteady Stokes system, $\bm{A} = \frac{1}{\alpha \Delta t}\bm{M}+\bm{L}$, so we can use the previously-defined Helmholtz preconditioner.
While these are vector linear systems, both $\bm{A}$ operators decouple the dimensions.
Therefore, applying the scalar preconditioners $d$ times (once in each dimension) directly provides a preconditioner for these vector linear systems.
Since the CG iterations do not correspond to a fixed linear operator, it is important that we make use of flexible GMRES as an outer iteration.
The sub-problems need not be solved exactly, and empirically two or three CG iterations are sufficient to provide an effective preconditioner.

For the steady Stokes system, $S = -D\bm{L}^{-1}\bm{G}$.
Before creating an approximate solver for $S$, we notice that $L = -\nu D\bm{M}^{-1} \bm{G}$.
To construct the approximate solver $\tilde{S}^{-1}$, we make the standard commutativity approximation $\bm L \bm M^{-1} \bm G \approx \bm G M^{-1} L$.
Then, we see that $\bm{L}^{-1} \bm{G} \approx \bm M^{-1} \bm{G} L^{-1} M$, and so
\begin{equation}
S = -D \bm{L}^{-1}\bm{G} \approx -D \bm{M}^{-1} \bm{G} L^{-1} M = \frac{1}{\nu} L L^{-1} M = \frac{1}{\nu} M.
\end{equation}
That is, the mass matrix $M$ provides an approximation to $S$.
Note that this approximation is, in fact, exact when the operators $\bm L$ and $\bm{G}$ commute, such as in the case of periodic boundary conditions.
The action of the approximate solver $\tilde{S}^{-1}$ is given by the diagonal mass preconditioner described in \hyperref[sec:massPrecond]{Section~\ref{sec:massPrecond}}.
In practice, approximating the action of $S^{-1}$ by the action of $\tilde{S}^{-1}$ doubles to triples the iterations of the iterative solver.
However, inverting $S$ is infeasible, whereas applying $\tilde{S}^{-1}$ is efficiently performed in $\mathcal{O}(p^d)$ time.

Likewise, for the unsteady Stokes system, the Schur complement is given by
\begin{equation}
S = -D \paren{\frac{1}{\alpha \Delta t}\bm{M}+\bm{L}}^{-1} \bm{G}.
\end{equation}
Using the same commutativity approximation as in the steady case, we obtain
\begin{equation}
\paren{\frac{1}{\alpha \Delta t}\bm{M}+\bm{L}}^{-1}\bm{G} \approx \bm{M}^{-1} \bm{G} \paren{\frac{1}{\alpha \Delta t}M+L}^{-1} M
\end{equation}
and so
\begin{equation}
S \approx -D\bm{M}^{-1}\bm{G} \paren{\frac{1}{\alpha \Delta t}M+L}^{-1} M = \frac{1}{\nu} L \paren{\frac{1}{\alpha \Delta t}M+L}^{-1} M.
\end{equation}
Therefore, we take
\begin{equation}
\tilde{S}^{-1} = \frac{\nu}{\alpha \Delta t} \tilde{L}^{-1} + \nu \tilde{M}^{-1}.
\label{eq:unsteadyStokesSchur}
\end{equation}
From \eqref{eq:unsteadyStokesSchur}, we can apply the action of $\tilde{S}^{-1}$ matrix-free by once again reusing the Poisson solver from \hyperref[sec:poissonPrecond]{Section~\ref{sec:poissonPrecond}} and the diagonal $\tilde{M}^{-1}$ from \hyperref[sec:massPrecond]{Section~\ref{sec:massPrecond}}.

\subsection{Relationship to other Schwarz-based solvers}
\label{sec:solverComparison}

A number of other matrix-free solvers based on a Schwarz methodology have been proposed for the solution of the high-order Poisson problem, and by extension the incompressible Navier-Stokes equations.
Closely related to the present work, multigrid solvers with matrix-free Schwarz-based smoothers for the spectral element method were constructed in \cite{Lottes2005}.
These methods were later extended in \cite{Fischer2005} to solve the unsteady Navier-Stokes equations, and were shown to perform efficiently on several large-scale tensor-product meshes.
The additive Schwarz smoothers used in \cite{Lottes2005} and \cite{Fischer2005} are constructed using tensor-product subdomains corresponding to the spectral elements of the high-order discretization.
Because each of the subdomains possesses a tensor-product geometry, the fast diagonalization method may be used to efficiently solve the local problems \cite{Wilhelm2001,Shen2011}.
In this work, however, we make use of subdomains defined by fully unstructured vertex patches, which in general do not possess a tensor-product structure.
The application of fast diagonalization methods to such geometries is non-trivial \cite{Pazner2018Approximate}.
For this reason, instead of using fast diagonalization to solve the local problems, we opt to solve the local low-order refined problems using an element-structured geometric multigrid V-cycle with ILU smoothing.

\section{GPU implementation}
\label{sec:gpu}

We have implemented the numerical algorithms described above in the framework of the MFEM finite element library~\cite{Anderson2019,MFEM}.
These algorithms take the form of single-source compute kernels that can target several different backends, including OpenMP on traditional CPUs as well as CUDA for use on the GPU.
In this section, we will describe important practical details that were required to obtain performant implementations for these algorithms.

\begin{figure*}[htbp!]
\hspace*{\fill}
\includegraphics{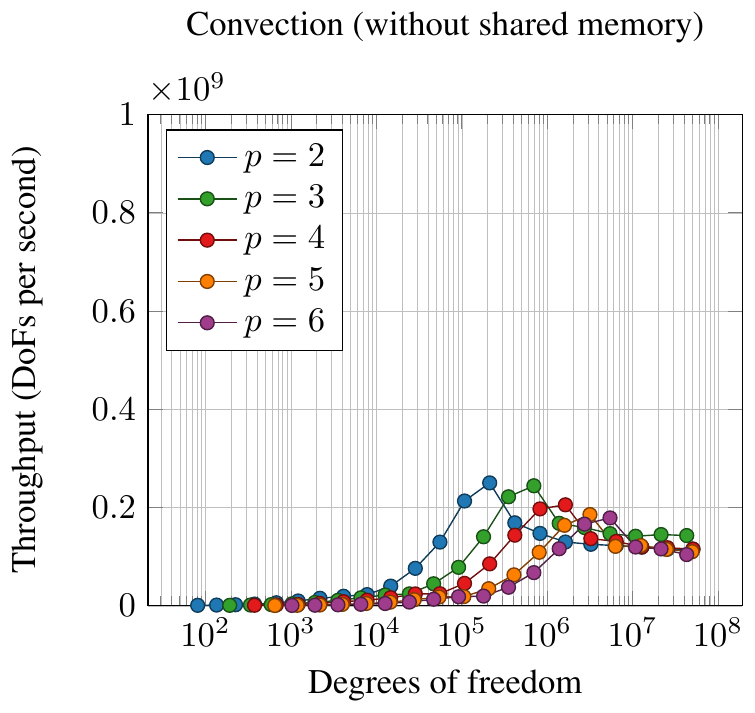}
\hspace*{\fill}
\includegraphics{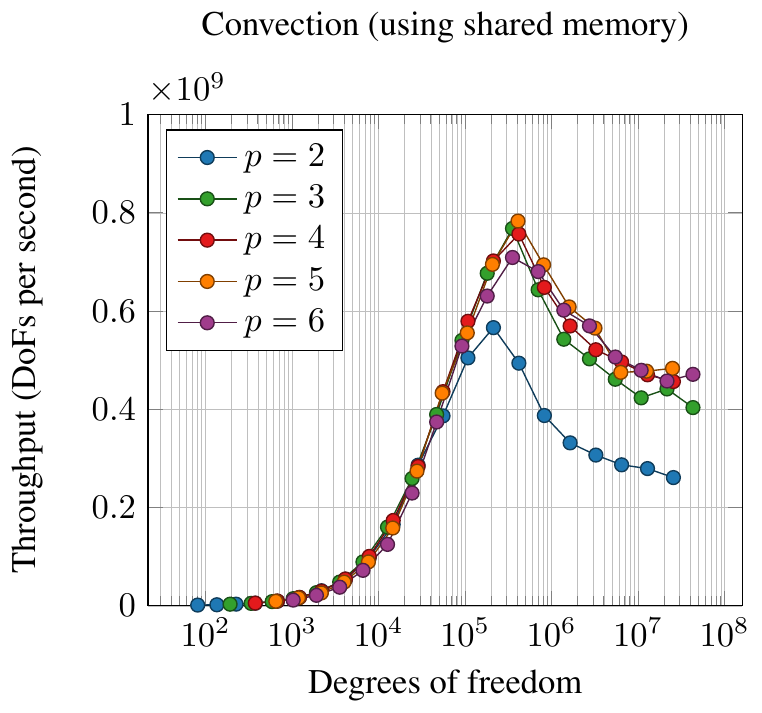}
\hspace*{\fill}

\caption{Throughput for nonlinear vector convection evaluation $\bm{N}(\bm{u})$ in 3D for several polynomial degrees on the GPU. Left: Initial implementation. Right: Elementwise shared-memory implementation. Efficiently reusing shared memory increases throughput by a factor of between 4 and 8 and also allows for the solution of larger problems.}
\label{fig:sharedVsNonShared}

\vspace{\floatsep}

\hspace*{\fill}
\includegraphics{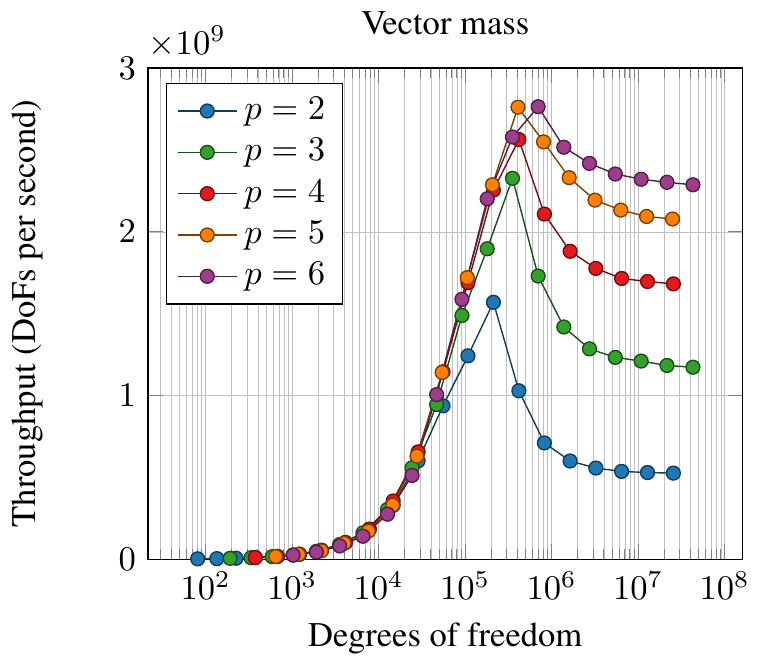}
\hspace*{\fill}
\includegraphics{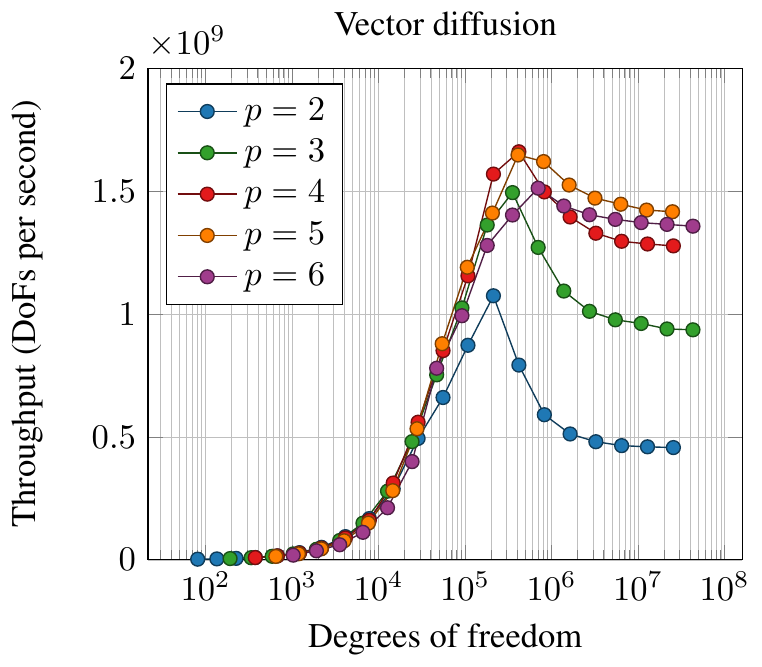}
\hspace*{\fill}

\vspace{\floatsep}

\hspace*{\fill}
\includegraphics{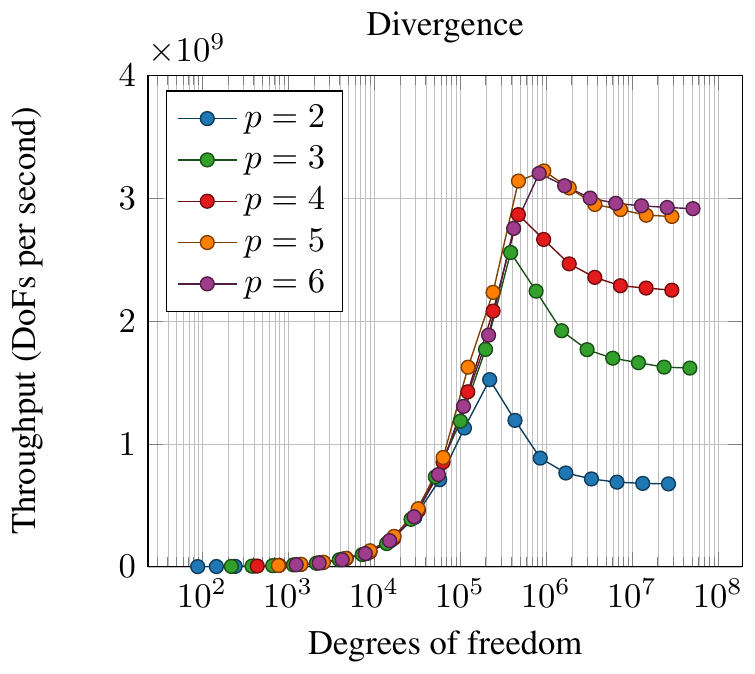}
\hspace*{\fill}
\includegraphics{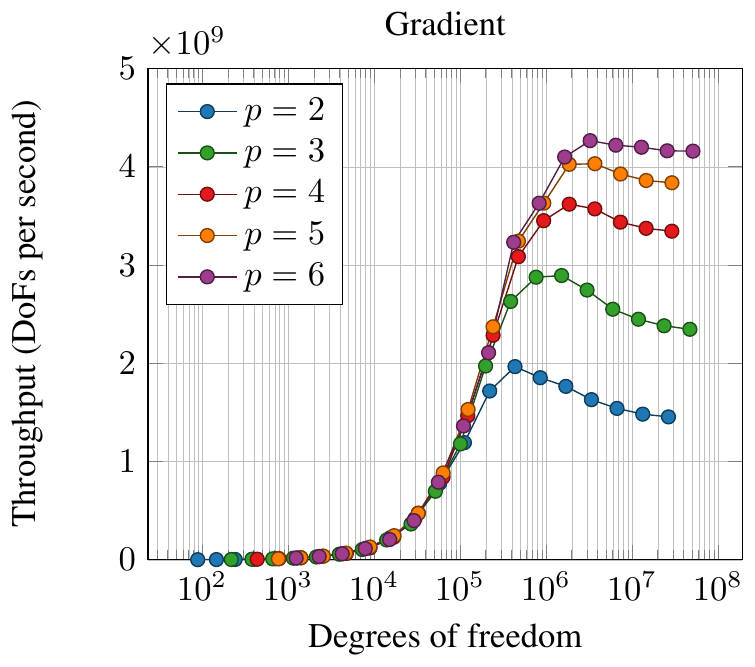}
\hspace*{\fill}
\caption{Throughput for linear operator evaluation kernels $\bm{M}$, $\bm{L}$, $\bm{D}$, $G$ in 3D for orders $p=2$ to $6$ on the GPU. Maximum throughput is achieved for higher orders ($p>3$) and larger problems (more than $10^6$ DoFs).}
\label{fig:matrixFreeRoofline}
\end{figure*}

Firstly, it is important to minimize memory movement between the CPU and GPU.
Modern GPU architectures have limited memory and cache sizes compared to their CPU counterparts, but can reach higher peak performance in terms of floating point operations.
This combination means that performance is only achievable if algorithms reach higher arithmetic intensities~\cite{Kronbichler2019}, thus motivating high-order matrix-free operators and solvers that have this potential.

We have found that using loop bounds known at compile time drastically improves performance.
In $p$-dependent compute kernels, these compile-time constants permit shared memory access within an element, thus reducing memory allocations and movement.
Figure~\ref{fig:sharedVsNonShared} shows the benefits of utilizing shared memory for the intra-element operations.
We see that the shared-memory implementation outperforms the naive implementation by a factor of $4$-$8$.
In practice, just-in-time compilation or explicit template instantiation can be used to allow for arbitrary-order simulations.

Figure~\ref{fig:matrixFreeRoofline} shows the throughput plots of our optimized implementations of the matrix-free linear operator evaluations in 3D.
These results were performed on a single Nvidia V100 GPU, showing the throughput achieved for various problem sizes and polynomial degrees.
Similar benchmark problems were used to assess kernel and backend performance in the context of the CEED project \cite{Dobrev2017}.

We believe that our implementation can be improved by better taking advantage of the GPU's shared memory.
First of all, symmetries exist in the 1D interpolation and differentiation matrices~\eqref{eq:B1D} and \eqref{eq:D1D} since the Gauss-Lobatto nodes are symmetric about the origin.
Taking advantage of these symmetries would approximately halve shared memory usage and simultaneously increase the arithmetic intensity of all operator evaluations.
In some sense, this improvement can be viewed as an extension of the original sum-factorization techniques, which use inherent operator symmetries to avoid extra memory storage.
Moreover, we see from Figure~\ref{fig:matrixFreeRoofline} that our implementation achieves lower throughput at lower orders.
It may be possible to address this issue by combining several low-order elements per thread block to yield higher throughput.
Unfortunately, the shared-memory requirements of our approach increase with order $p$ and dimension $d$, so exhausting the shared memory is inevitable with very high order simulations.
If necessary, one could compensate by storing $B$ and $D$ in the GPU's global memory to allow for even higher order simulations.

\section{Numerical results and performance analysis}
\label{sec:numericalResults}

\subsection{Sub-problem solver performance}
\label{sec:subproblemResults}

We first assess the performance of the matrix-free sub-problem preconditioners by measuring the number of Krylov iterations required to converge to a fixed tolerance under $h$- and $p$-refine\-ment.
We solve the linear system $A\bm x = \bm b$, where $A$ is either the mass matrix, Laplacian operator, or positive-definite Helmholtz operator.
We write the Helmholtz operator as $M/\Delta t + L$, and choose two representative time steps: $\Delta t = 10^{-1}$ and $\Delta t = 10^{-3}$.
For each of the problems considered, we use a preconditioned conjugate gradient iteration to solve the problem, with a relative residual tolerance of $10^{-8}$ as a stopping criterion.
The right hand side $\bm b$ is taken to be a random vector.
For the mass matrix, the preconditioner is the collocated diagonal preconditioner described above, and for the Helmholtz and Poisson solvers, we use the low-order refined parallel subspace correction procedure.

\begin{figure*}
\hspace*{\fill}
\includegraphics{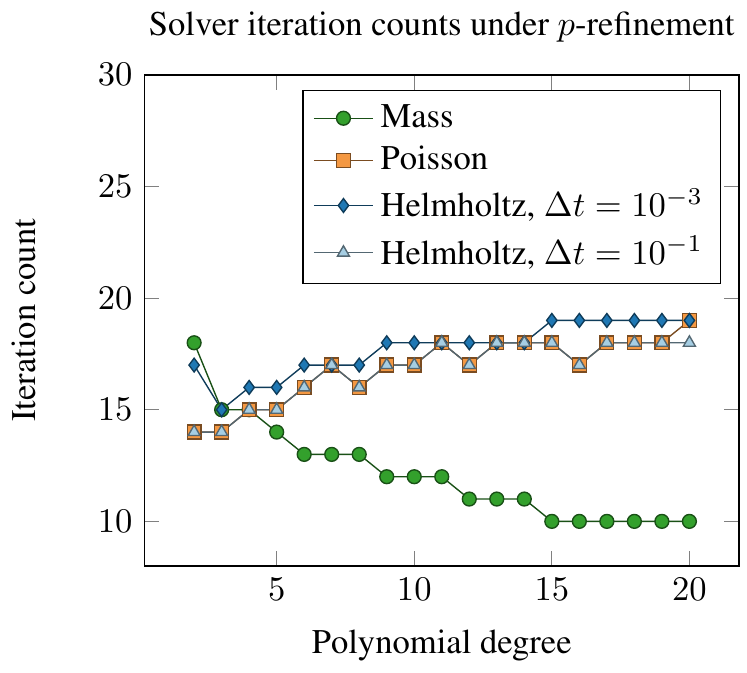}
\hspace*{\fill}
\includegraphics{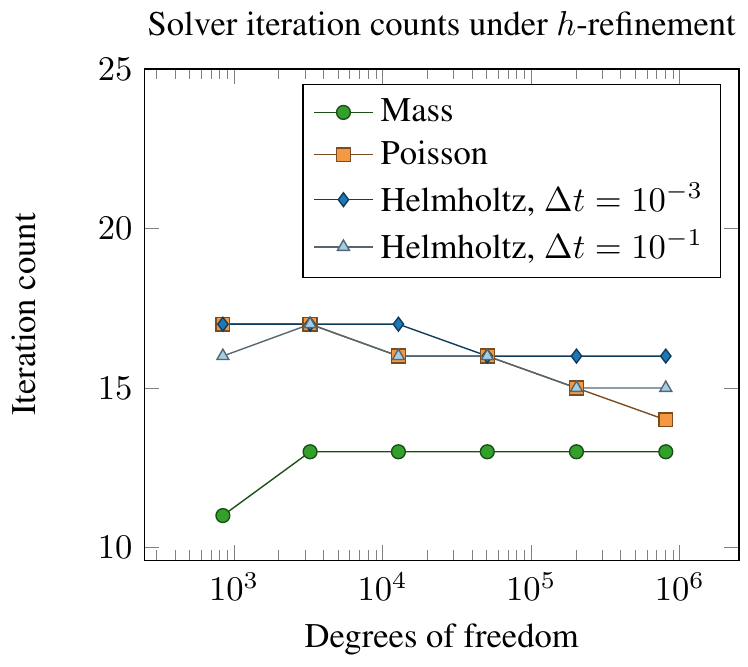}
\hspace*{\fill}
\caption{Iteration counts for sub-problem solvers under $p$- and $h$-refinement.
For the case of $h$-refinement, we use a fixed polynomial degree of $p=7$.}
\label{fig:sub-problem-iters}
\end{figure*}

To perform the $p$-refinement study, we use a fixed Cartesian grid with 64 elements in two dimensions, and polynomial degrees from $p=2$ to $p=20$.
The number of iterations required to converge to tolerance is shown in Figure \ref{fig:sub-problem-iters}.
We note that the number of iterations remains bounded for all problems and for all polynomial degrees.
We observe a slight pre-asymptotic increase in the number of iterations for the Poisson and Helmholtz problems, but the iteration counts remain below 20 for all cases.
The number of iterations required for the mass solve decreases with increasing polynomial degree.
This is corroborated by an eigenvalue analysis, which shows decreasing condition number of the preconditioned system with increasing polynomial degree.

For the case of $h$-refinement, we fix the polynomial degree to be $p=7$, and perform a sequence of uniform refinements.
The initial mesh is a $4 \times 4$ Cartesian grid with 841 DoFs.
We perform five refinements, so that the finest mesh is $128 \times 128$ with 804{,}609 DoFs.
The number of iterations required to converge to tolerance is shown in Figure \ref{fig:sub-problem-iters}.
Here, we observe approximately constant iterations, independent of the mesh refinement.
These examples verify the robustness of the sub-problem preconditioners with respect to the mesh size and polynomial degree.

\subsection{Steady-state Stokes flow}
\label{sec:steadyStokesResults}

\begin{figure}
\includegraphics{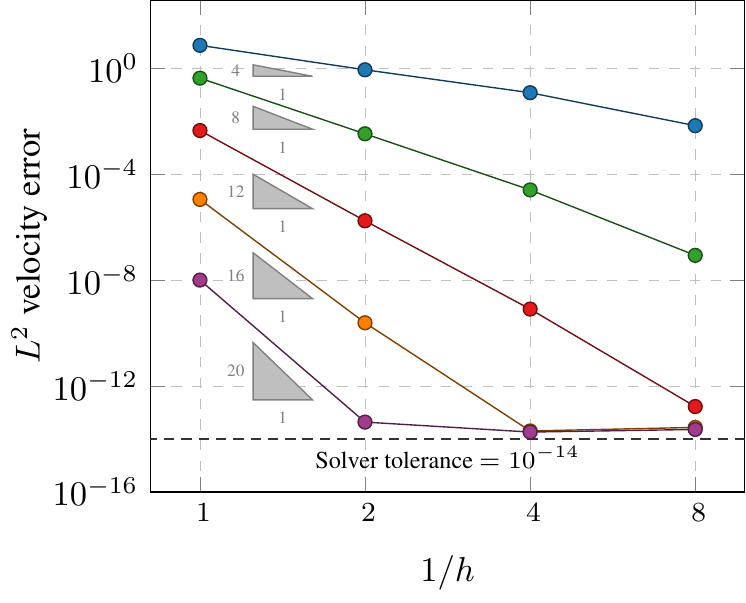}
\caption{$L^2$ velocity error showing high-order spatial convergence for steady-state Stokes. Polynomial degrees $3, 7, 11, 15,$ and $19$ are used for the velocity finite element space.}
\label{fig:stokes-convergence}
\end{figure}

To verify the high-order accuracy of the spatial discretization, and to test the convergence properties of the solver, we solve the steady Stokes equations with a smooth solution in two spatial dimensions.
Setting $\nu = 1$, we choose the right-hand side to be
\begin{align*}
   f_1(x,y) &= \pi  \cos (\pi  y) \left(4 \pi ^2 (1-2 \cos (2 \pi  x)) \sin (\pi  y)-\sin (\pi
   x)\right), \\
   f_2(x,y) &= 2 \pi ^3 \sin (2 \pi  x) (2 \cos (2 \pi  y)-1)-\pi  \cos (\pi  x) \sin (\pi  y).
\end{align*}
The exact solution is then given by
\begin{align*}
u_1(x,y) &= 2 \pi  \sin ^2(\pi  x) \sin (\pi  y) \cos (\pi  y), \\
u_2(x,y) &= -2 \pi  \sin (\pi  x) \cos (\pi  x) \sin ^2(\pi  y), \\
p(x,y) &= \cos (\pi  x) \cos (\pi  y).
\end{align*}
The exact solution is imposed as a Dirichlet boundary condition for velocity on all domain boundaries.
We run this case with polynomial degrees from $p=3$ to $p=19$ on a sequence of uniformly refined Cartesian meshes.
We solve the resulting linear system using the FGMRES method with the matrix-free block triangular preconditioners described in \hyperref[sec:blockPrecond]{Section~\ref{sec:blockPrecond}}.
The stopping criterion for the iterative solver is a relative residual norm of $10^{-14}$.
The spatial convergence is shown in Figure \ref{fig:stokes-convergence}.
The expected $p+1$ order of accuracy was observed in all cases, verifying high-order spatial accuracy.
In Table \ref{tab:stokes-mms}, we show the $L^2$ error for velocity and pressure for cases considered, together with the number of FMGRES iterations required to converge the solution.
The number of iterations shows a slight pre-asymptotic increase, but remains bounded with respect to both $h$ and $p$.

\begin{table*}
\caption{Error and convergence results for steady Stokes equation, showing $L^2$ error norms for velocity and pressure, and number of FGMRES iterations required to reduce the residual by a factor of $10^{14}$.}
\label{tab:stokes-mms}

\hspace*{\fill}
\begin{tabular}{rlllll}
\multicolumn{6}{c}{\parbox[t][4ex]{1cm}{$p=7$}}\\
\toprule
$1/h$ & $\|\bm u_h - \bm u\|_2$ & Rate & $ \|p_h - p \|_2$ & Rate & Its. \\
\midrule
1 & $4.25\times10^{-1}$ & --- & $1.61\times10^{-1}$ & --- & 31 \\
4 & $3.36\times10^{-3}$ & 6.99 & $3.33\times10^{-3}$ & 5.60 & 41 \\
16 & $2.55\times10^{-5}$ & 7.04 & $7.58\times10^{-5}$ & 5.46 & 43 \\
64 & $8.63\times10^{-8}$ & 8.21 & $2.41\times10^{-7}$ & 8.30 & 46 \\
\bottomrule
\end{tabular}
\hspace*{\fill}
\begin{tabular}{rlllll}
\multicolumn{6}{c}{\parbox[t][4ex]{1cm}{$p=11$}}\\
\toprule
$1/h$ & $\|\bm u_h - \bm u\|_2$ & Rate & $ \|p_h - p \|_2$ & Rate & Its. \\
\midrule
1 & $4.46\times10^{-3}$ & --- & $1.72\times10^{-2}$ & --- & 35 \\
4 & $1.74\times10^{-6}$ & 11.32 & $1.33\times10^{-6}$ & 13.66 & 42 \\
16 & $8.05\times10^{-10}$ & 11.08 & $2.57\times10^{-9}$ & 9.01 & 47 \\
64 & $1.68\times10^{-13}$ & 12.23 & $6.76\times10^{-13}$ & 11.89 & 49 \\
\bottomrule
\end{tabular}
\hspace*{\fill}

\vspace{\floatsep}

\hspace*{\fill}
\begin{tabular}{rlllll}
\multicolumn{6}{c}{\parbox[t][4ex]{1cm}{$p=15$}}\\
\toprule
$1/h$ & $\|\bm u_h - \bm u\|_2$ & Rate & $ \|p_h - p \|_2$ & Rate & Its. \\
\midrule
1 & $1.12\times10^{-5}$ & --- & $1.05\times10^{-2}$ & --- & 39 \\
4 & $2.44\times10^{-10}$ & 15.48 & $1.46\times10^{-10}$ & 26.10 & 45 \\
16 & $2.01\times10^{-14}$ & 13.57 & $5.30\times10^{-13}$ & 8.11 & 52 \\
64 & $2.73\times10^{-14}$ & -0.44 & $7.43\times10^{-13}$ & -0.49 & 53 \\
\bottomrule
\end{tabular}
\hspace*{\fill}
\begin{tabular}{rlllll}
\multicolumn{6}{c}{\parbox[t][4ex]{1cm}{$p=19$}}\\
\toprule
$1/h$ & $\|\bm u_h - \bm u\|_2$ & Rate & $ \|p_h - p \|_2$ & Rate & Its. \\
\midrule
1 & $1.00\times10^{-8}$ & --- & $7.24\times10^{-3}$ & --- & 41 \\
4 & $4.34\times10^{-14}$ & 17.82 & $3.22\times10^{-13}$ & 34.39 & 47 \\
16 & $1.82\times10^{-14}$ & 1.25 & $7.48\times10^{-13}$ & -1.21 & 52 \\
64 & $2.29\times10^{-14}$ & -0.33 & $8.69\times10^{-13}$ & -0.22 & 55 \\
\bottomrule
\end{tabular}
\hspace*{\fill}
\end{table*}

\subsection{Unsteady Stokes flow}
\label{sec:unsteadyStokesResults}

\begin{figure*}
\centering
\hspace*{\fill}
\includegraphics{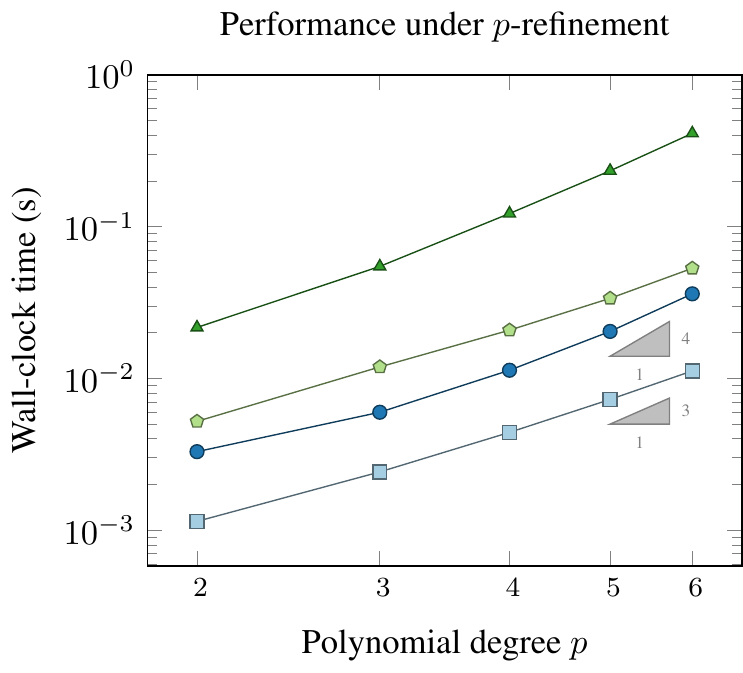}
\hspace*{\fill}
\includegraphics{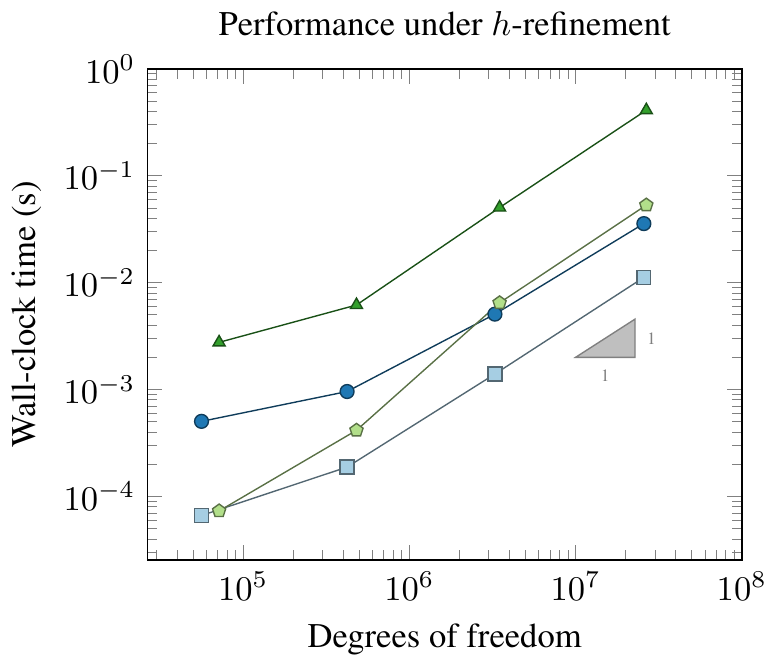}
\hspace*{\fill}

\includegraphics{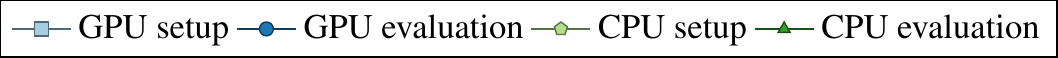}
\caption{Performance of unsteady Stokes operator in 3D. Our implementation achieves expected rates of $\mathcal{O}(p^d)$ and $\mathcal{O}(p^{d+1})$ for matrix-free setup and evaluation of the block operator, respectively. Shared-memory GPU implementation of operator evaluation outperforms 20-core CPU implementation by a factor of 6 at $p=2$ and 11 at $p=6$.}
\label{fig:stokes-timings}
\end{figure*}

We evaluate the performance of our GPU implementation of the coupled unsteady Stokes operator defined by~\eqref{eq:unsteadyStokesDiscrete}.
Figure~\ref{fig:stokes-timings} records the 3D operator setup and evaluation times under $p$-refinement for a fixed mesh with 32{,}768 hexahedral elements.
The GPU results were computed using one Nvidia V100 GPU, while the CPU results used 20 POWER8 cores on one node of Lawrence Livermore National Laboratory's Ray supercomputer.
Our implementation achieves the expected rates of $\mathcal{O}(p^d)$ and $\mathcal{O}(p^{d+1})$ for matrix-free setup and evaluation, respectively.
Moreover, we can see the computational benefits of high-order simulations on the GPU, as our optimized GPU kernels outperform the 20-core CPU implementation by a factor of $6$ at $p=2$ and $11$ at $p=6$.

Figure~\ref{fig:stokes-timings} also shows the performance under $h$-refinement, fixing the polynomial degree at $p=6$.
We see that for fixed $p$, the wall-clock time scales linearly with DoFs.
Once the device is saturated, the operator evaluation is accelerated by a factor of about 11 relative to the 20-core CPU evaluation.

\subsection{Incompressible Navier-Stokes: Kovasznay flow}
\label{sec:kovasznay}

\begin{figure}[t]
\begin{minipage}[t]{0.45\linewidth}
\centering
\includegraphics[width=0.55\linewidth]{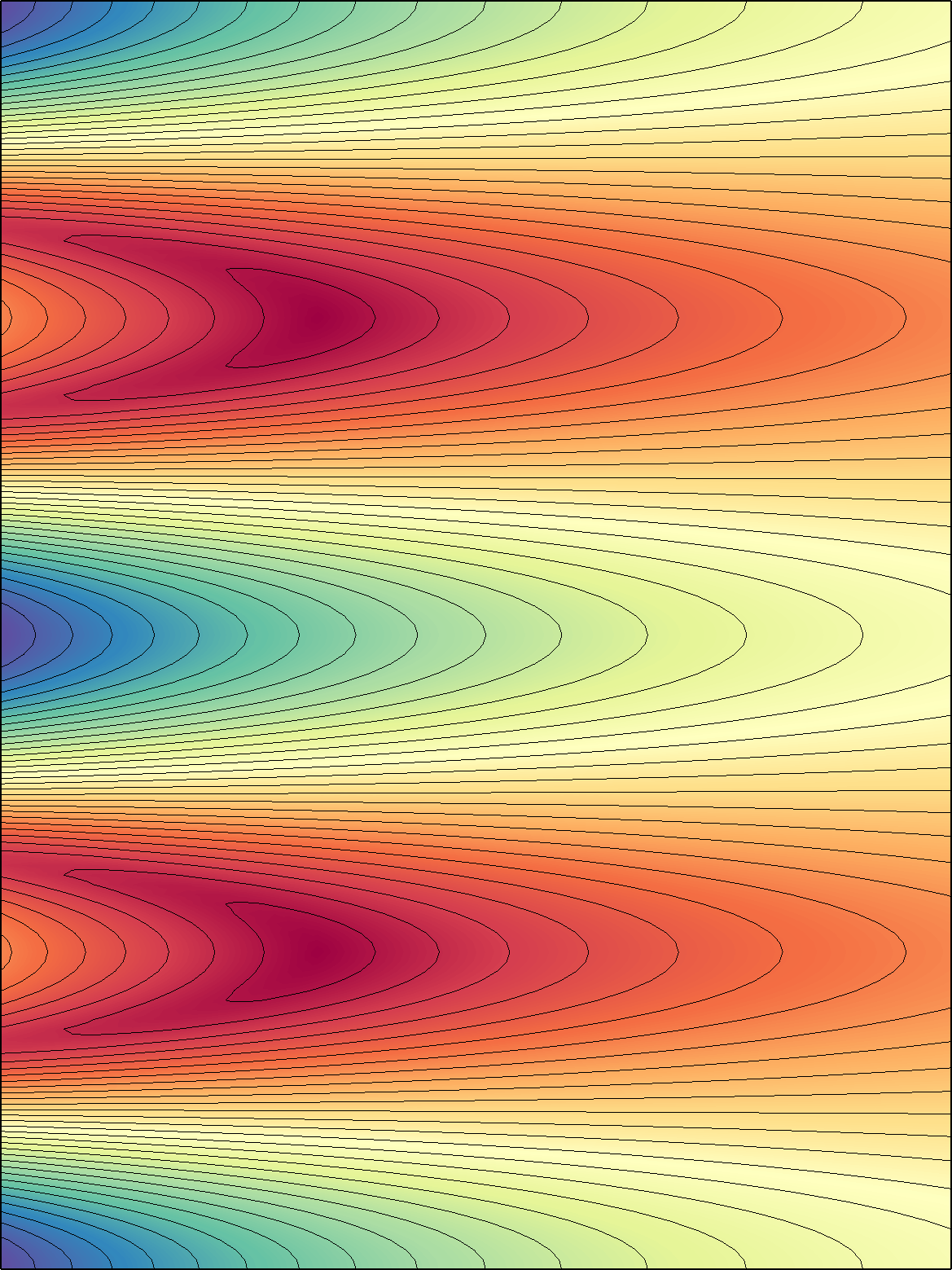}
\caption{Two-dimensional Kovasznay flow, showing contours of velocity magnitude computed using a coarse mesh with degree 11 polynomials.}
\label{fig:kov-soln}
\end{minipage}
\begin{minipage}[t]{0.45\linewidth}
\includegraphics{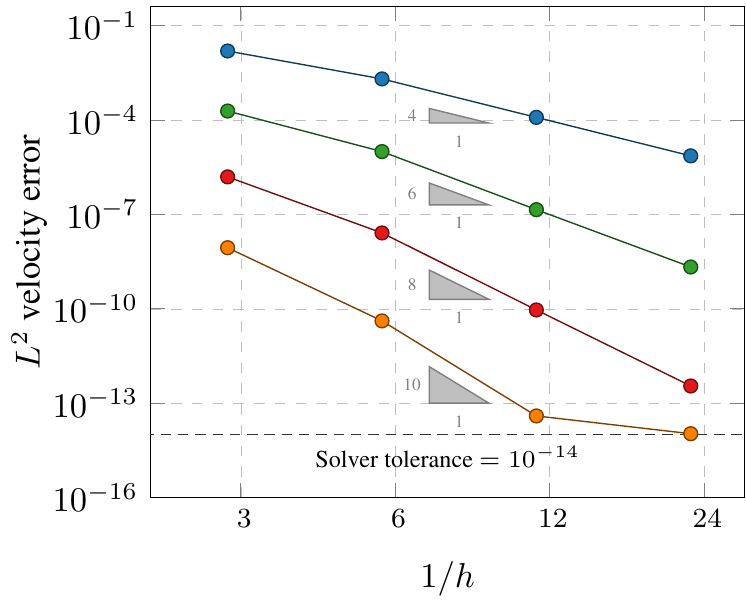}
\caption{Two-dimensional Kovasznay flow. $L^2$ velocity error using polynomial degrees $p=3,5,7,9$.}
\label{fig:kov}
\end{minipage}
\end{figure}

An analytical solution to the stationary Navier-Stokes equations in two spatial dimensions due to Kovasznay can be found in~\cite{Kovasznay1948}.
This solution may be used to represent the wake behind a periodic array of cylinders in the $y$-direction.
The solution is given by
\[\begin{aligned}
    \lambda &= \frac{\Re}{2} - \sqrt{\frac{\Re^2}{4} + 4 \pi^2},\\
    u_1(x,y) &= 1 - \exp{(\lambda x)} \cos{(2 \pi y)},\\
    u_2(x,y) &= \frac{\lambda}{2 \pi} \sin{(2 \pi x)}, \\
    p(x,y) &= -\frac{1}{2} \exp(2 \lambda x ),
\end{aligned}\]
where $\Re$ is the Reynolds number for the flow. In our example we define the problem in the rectangular domain $[-1/2, 1] \times [-1/2, 3/2]$ with $\Re = 40$.
Velocity magnitude contours of the solution are shown in Figure \ref{fig:kov-soln}.
We use pseudo time integration in order to apply the projection method described in~\hyperref[sec:projection]{Section~\ref{sec:projection}} to a steady-state problem.
The pseudo time step is chosen to be $\Delta t = 10^{-3}$ on the coarsest mesh, and is reduced by a factor of two with each refinement.
The exact solution is enforced as a Dirichlet boundary condition on all boundaries of the domain.
The equations are integrated until a final time of $t=8$ to allow for the errors to propagate out of the domain.

To investigate spatial convergence we compute the solution with increasing polynomial degree and uniform spatial refinement.
The $L^2$ errors for the velocity are shown in Figure~\ref{fig:kov}.
We observe the desired order of accuracy for all of the cases considered.

\subsection{Incompressible Navier-Stokes: Taylor-Green vortex}
\label{sec:taylorGreen}

We next consider the incompressible Taylor-Green vortex, which is a standard benchmark case often used to assess the accuracy of high-order methods \cite{Brachet1983,Brachet1984}.
This problem represents a simple model for the development of turbulence and resulting cascade of energy from large to small scales \cite{Shu2005}.
We use the problem configuration as defined in the first international workshop on high-order CFD methods \cite{Wang2013}.
The domain is taken to be the fully periodic cube $[-\pi, \pi]^3$ and the initial state is set to
\[\begin{aligned}
    u_1(x,y,z) &= \sin{(x)} \cos{(y)} \sin{(z)},\\
    u_2(x,y,z) &= -\cos{(x)} \sin{(y)} \sin{(z)},\\
    u_3(x,y,z) &= 0.
\end{aligned}\]
The Reynolds number is chosen to be $\Re = 1600$.
The equations are integrated with BDF-3 until a final time of $t=20$ using a time step of $\Delta t = 2.5 \times 10^{-3}$.
We consider the following three configurations:
\begin{itemize}
\item $p=3$, $24 \times 24 \times 24$ grid, 439{,}276 DoFs per component.
\item $p=7$, $12 \times 12 \times 12$ grid, 650{,}701 DoFs per component.
\item $p=11$, $6 \times 6 \times 6$ grid, 346{,}969 DoFs per component.
\end{itemize}

\begin{figure*}
\def\tgvwidth{0.3\linewidth}
\hspace*{\fill}
\begin{minipage}{\tgvwidth}
\centering
\includegraphics[width=\linewidth]{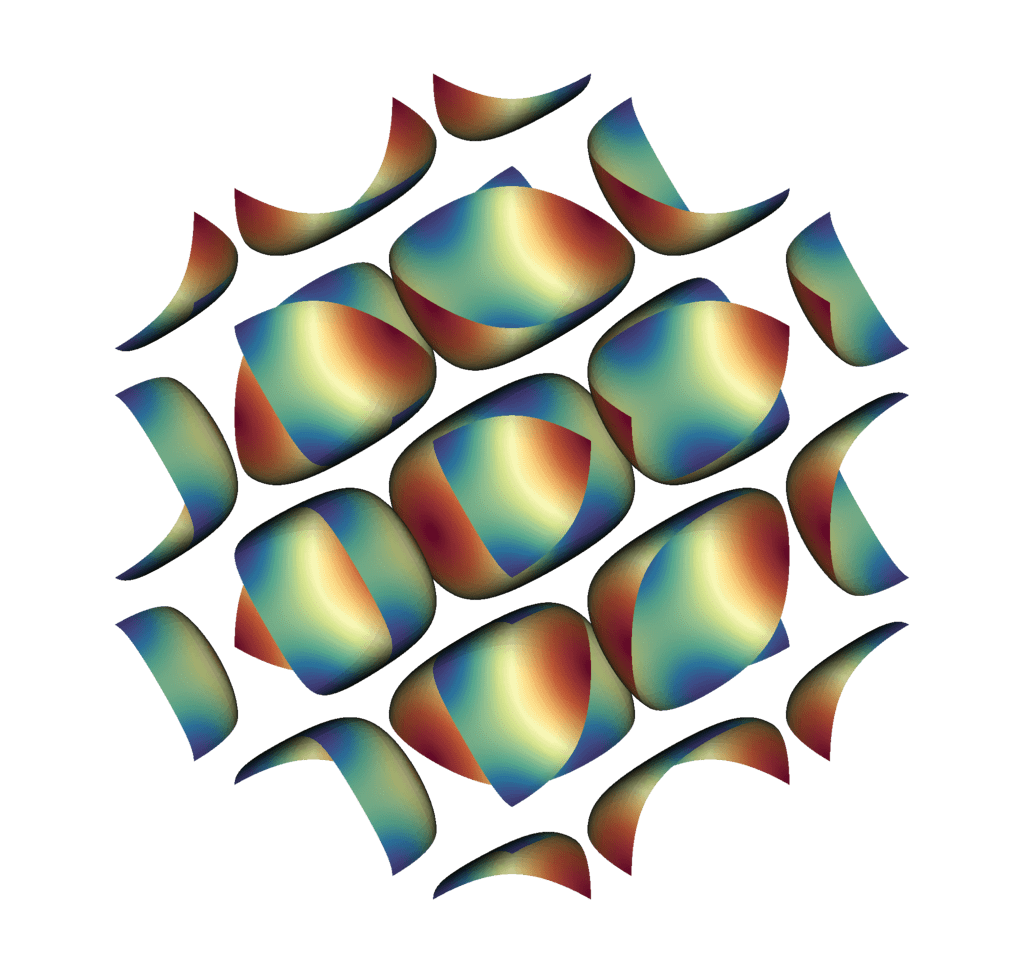} \\
$t=0$
\end{minipage}
\hspace*{\fill}
\begin{minipage}{\tgvwidth}
\centering
\includegraphics[width=\linewidth]{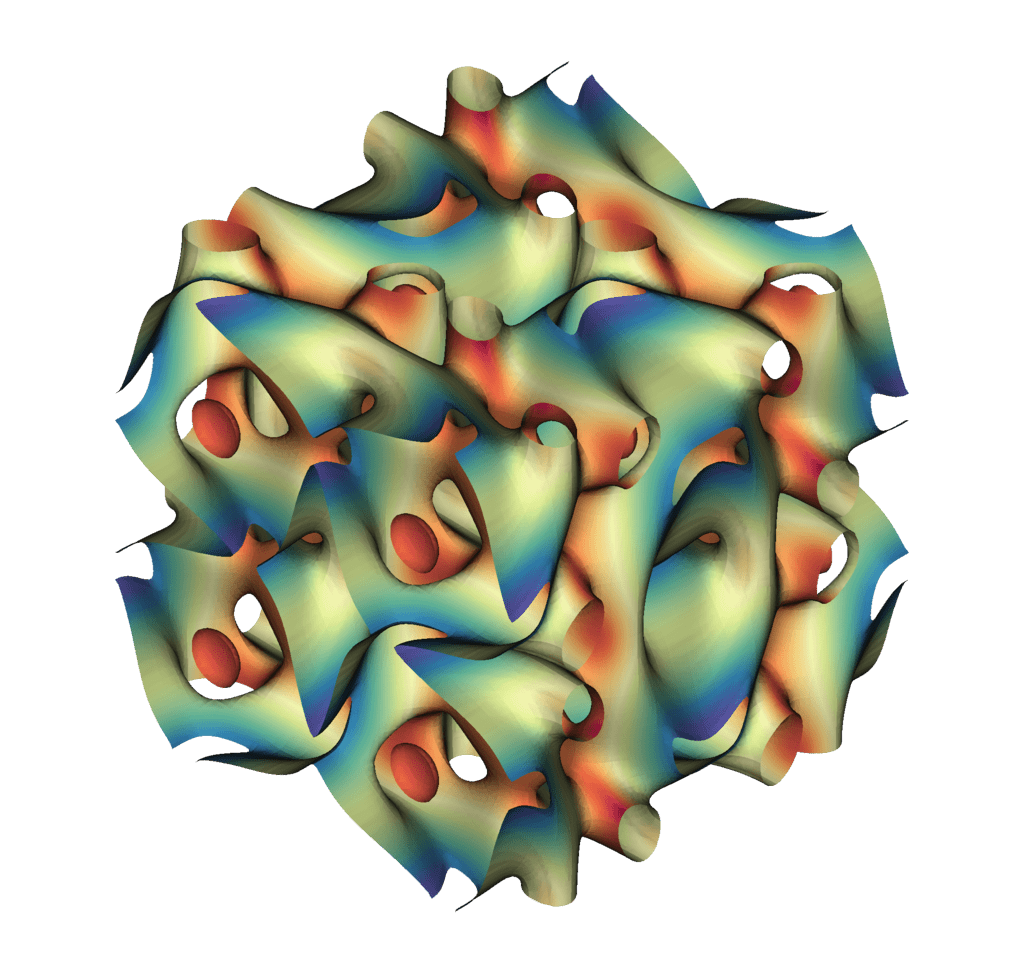} \\
$t=2$
\end{minipage}
\hspace*{\fill}

\vspace{\floatsep}

\hspace*{\fill}
\begin{minipage}{\tgvwidth}
\centering
\includegraphics[width=\linewidth]{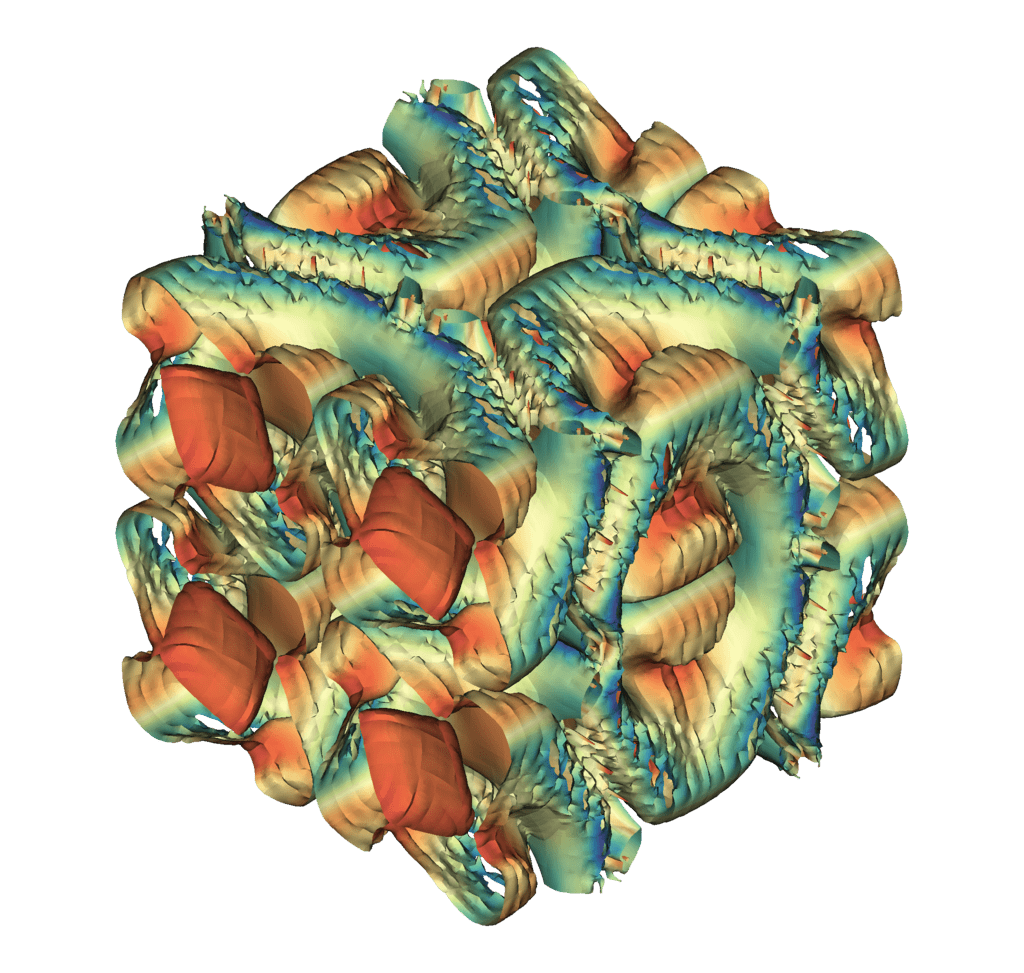} \\
$t=4$
\end{minipage}
\hspace*{\fill}
\begin{minipage}{\tgvwidth}
\centering
\includegraphics[width=\linewidth]{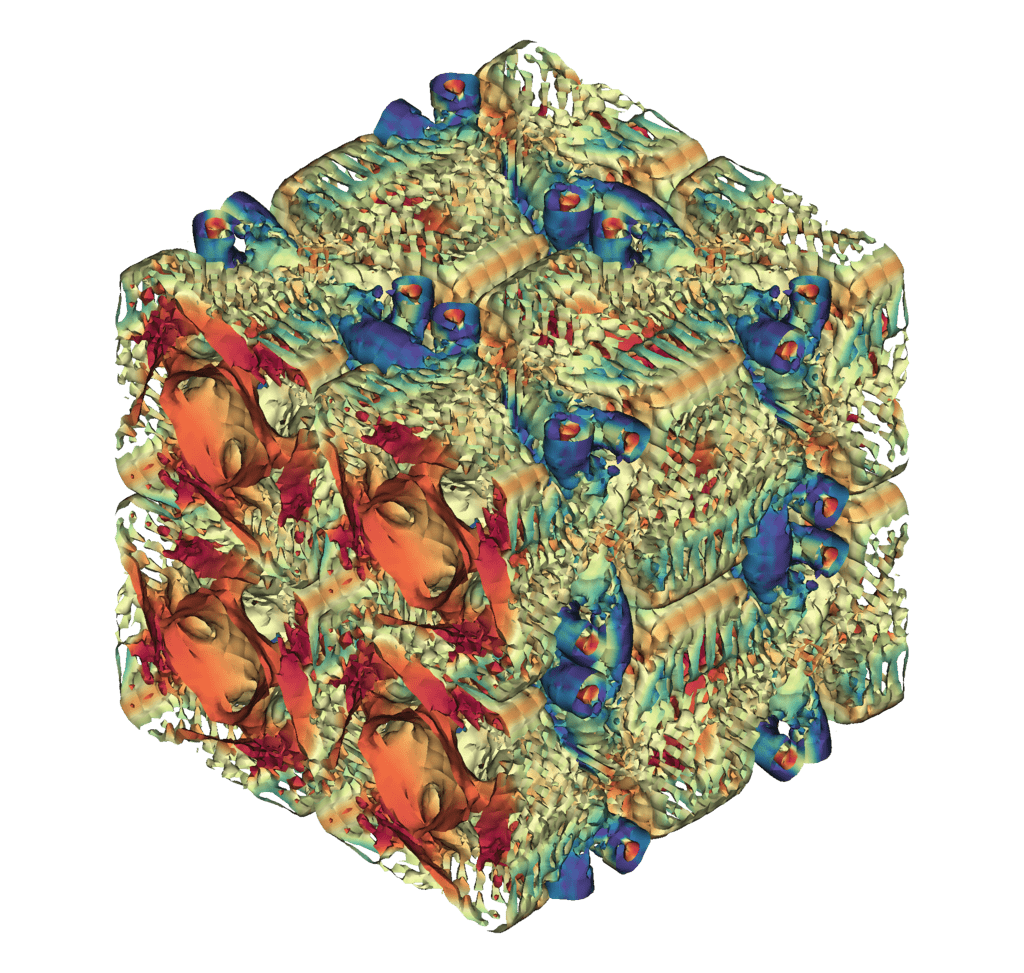} \\
$t=6$
\end{minipage}
\hspace*{\fill}

\vspace{\floatsep}

\hspace*{\fill}
\begin{minipage}{\tgvwidth}
\centering
\includegraphics[width=\linewidth]{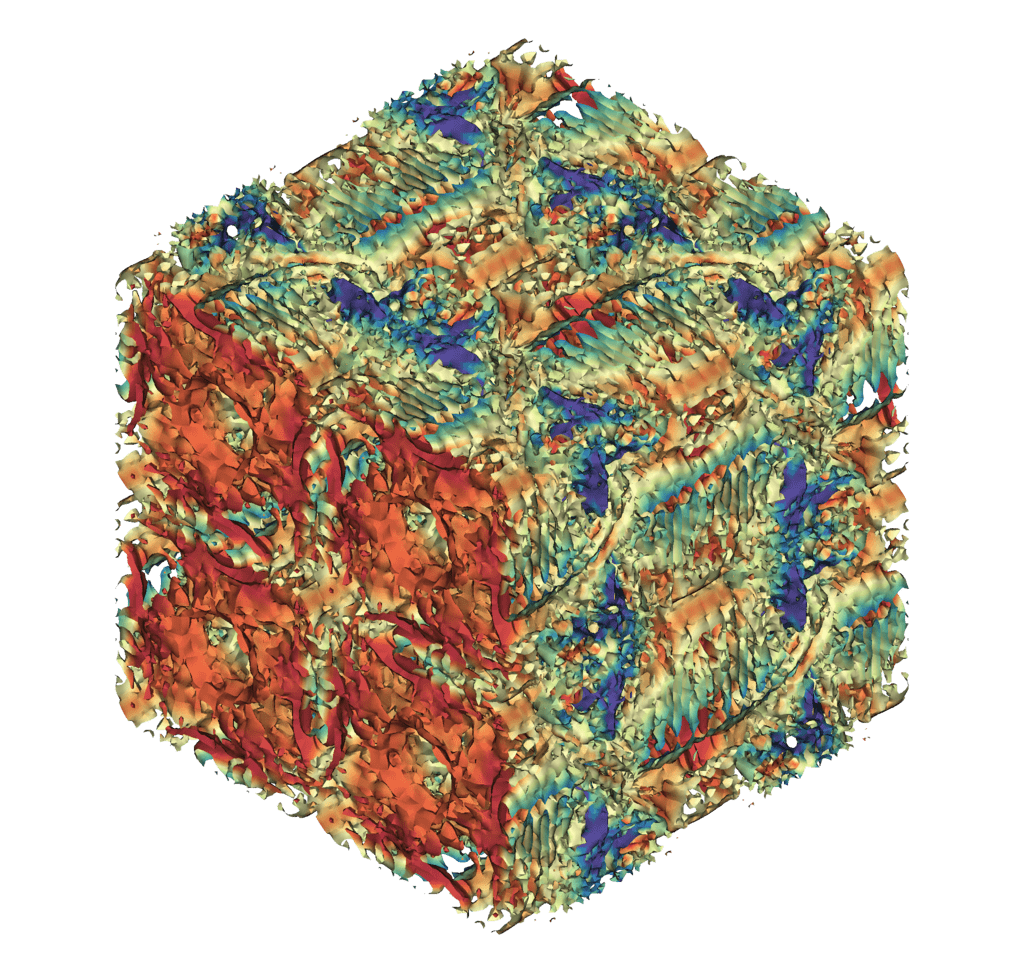} \\
$t=8$
\end{minipage}
\hspace*{\fill}
\begin{minipage}{\tgvwidth}
\centering
\includegraphics[width=\linewidth]{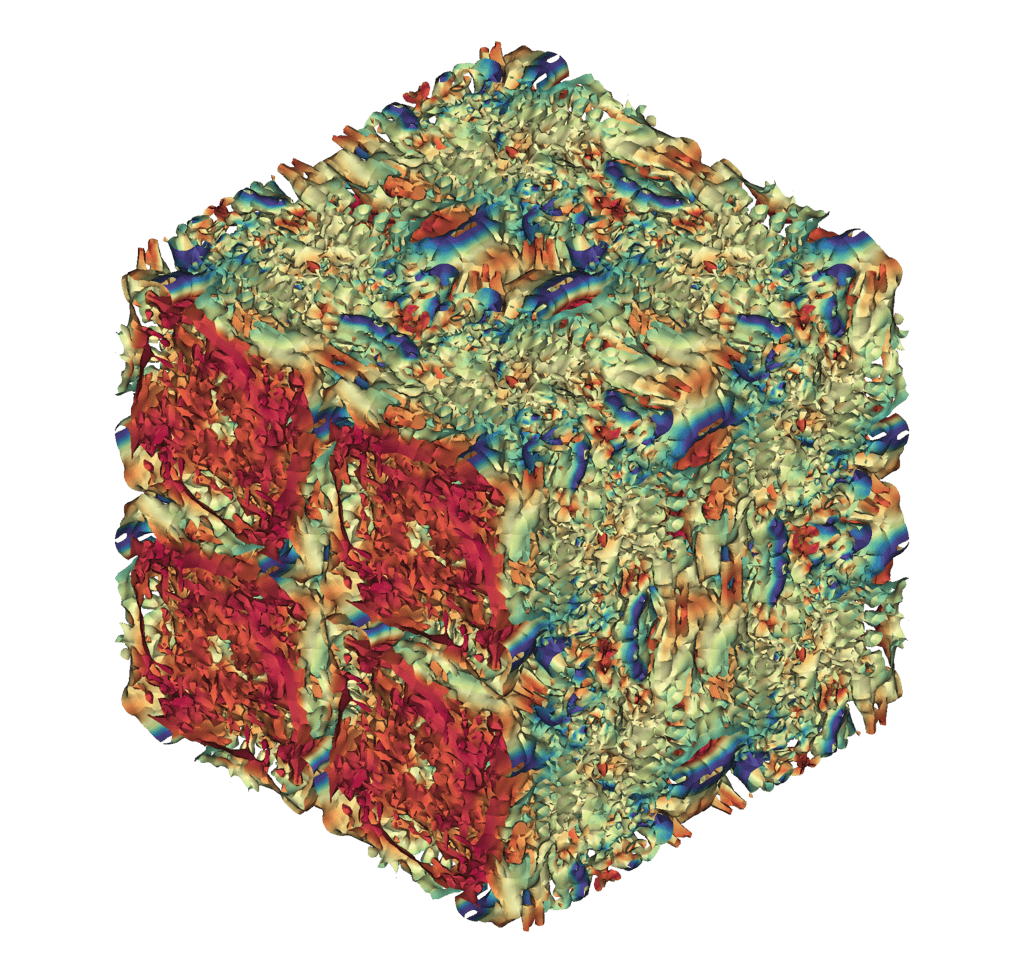} \\
$t=10$
\end{minipage}
\hspace*{\fill}

\caption{Time evolution of the incompressible Taylor-Green vortex, showing $Q=0.1$ isosurfaces colored by velocity magnitude.}
\label{fig:tgv-iso}
\end{figure*}

In Figure \ref{fig:tgv-iso}, we display the time evolution of $Q$-criterion isosurfaces, colored by velocity magnitude.
The quantity $Q$ is defined by $Q = \frac{1}{2} \sum_{i,j=1}^3 \partial_i u_j \partial_j u_i $ and is commonly used for vortex identification \cite{Hunt1988, Dubief2000}.
These isosurfaces clearly display the evolution from smooth, large-scale structures to small-scale turbulent structures.

\begin{figure}
\includegraphics{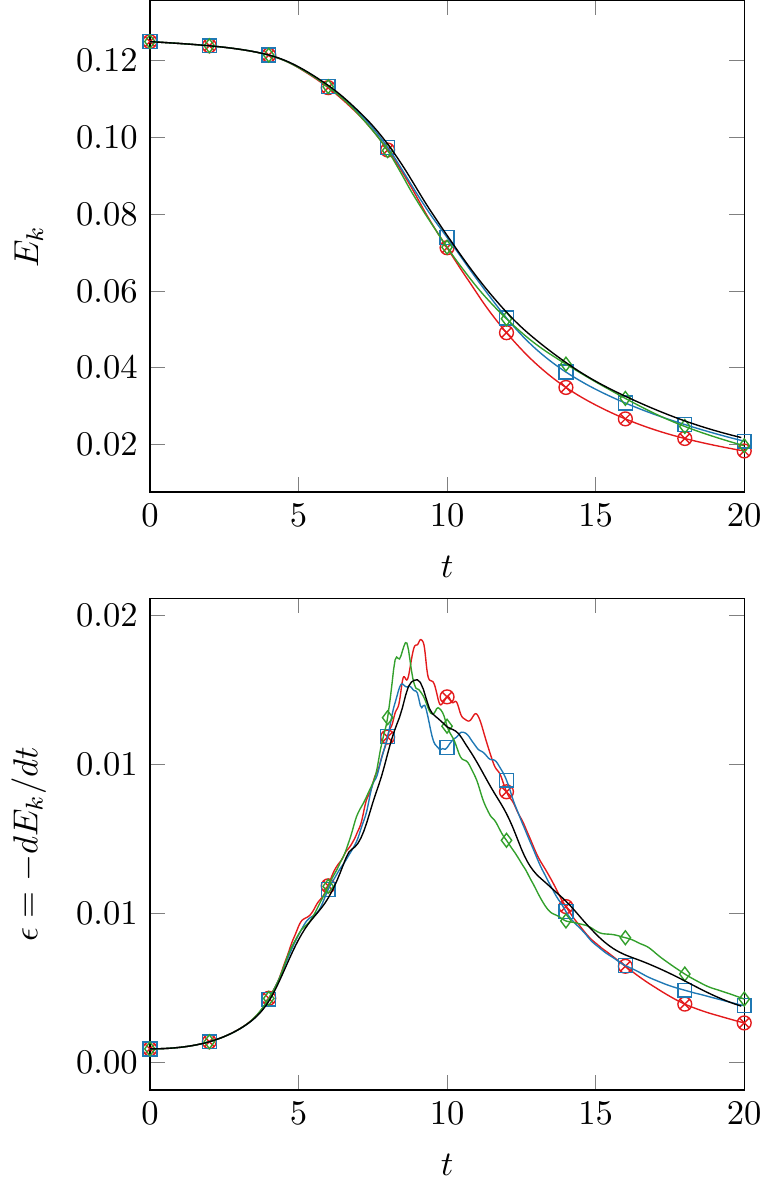}

\begin{center}
\includegraphics{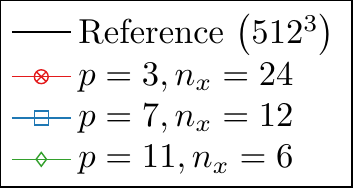}
\end{center}
\caption{Time evolution of total kinetic energy and kinetic energy dissipation rate for the incompressible Taylor-Green vortex. Comparison with reference data from a fully-resolved pseudo-spectral method with $512^3$ degrees of freedom per component.}
\label{fig:tgv-qoi}
\end{figure}

We compare the results obtained using the present high-order finite element flow solver with reference data obtained using a de-aliased pseudo-spectral method \cite{vanRees2011}.
The quantities of interest for this comparison are the total kinetic energy
\begin{align}
    E_k = \frac{1}{|\Omega|}\int_{\Omega} \frac{\bm u \cdot \bm u}{2} d\bm x
\end{align}
and the kinetic energy dissipation rate
\begin{align}
    \epsilon = -\frac{d E_k}{dt}.
\end{align}
The time evolution of these quantities is shown in Figure \ref{fig:tgv-qoi}.
All of the configurations considered result in good agreement with the reference data.
The lowest order case (polynomial degree $p=3$) is the most dissipative, slightly under-predicting the total kinetic energy after about $t=10$.
The highest order case we considered (polynomial degree $p=11$) gives results of comparable accuracy to the $p=7$ case, with roughly half as many degrees of freedom.

\section{Conclusions}
\label{sec:conclusions}

In this work, we have described a high-order finite element method for solving incompressible flow problems with a matrix-free implementation that efficiently utilizes the high performance and memory bandwidth of modern graphics processing units.
The resulting finite element operators are applied using sum-factorization algorithms and optimized GPU kernels, implemen\-ted in the MFEM library \cite{MFEM}, and obtain throughput of over a billion degrees of freedom per second on a single Nvidia V100 GPU.
A suite of matrix-free linear solvers was developed for the resulting mass, Poisson, Helmholtz, and Stokes linear systems.
These preconditioners do not require the costly assembly of high-order system matrices.
Their memory usage is optimal, requiring only constant memory per degree of freedom.
Furthermore, the number of operations required to apply the preconditioners scales linearly with the number of degrees of freedom, which is the same as the operator evaluation.
The robustness of preconditioners with respect to the mesh size, polynomial degree, and time step was demonstrated on a range of test problems.
The flow solver was applied to a variety of 2D and 3D problems, verifying the high-order accuracy and efficiency of the method.

\section*{Acknowledgements}
This material is based upon work supported by the National Science Foundation Graduate Research Fellowship Program under Grant No. DGE 1752814.
Any opinions, findings, and conclusions or recommendations expressed in this material are those of the authors and do not necessarily reflect the views of the National Science Foundation.
This work was performed under the auspices of the U.S. Department of Energy by Lawrence Livermore National Laboratory under Contract DE-AC52-07NA27344 (LLNL-JRNL-791975).

This document was prepared as an account of work sponsored by an agency of the United States government.
Neither the United States government nor Lawrence Livermore National Security, LLC, nor any of their employees makes any warranty, expressed or implied, or assumes any legal liability or responsibility for the accuracy, completeness, or usefulness of any information, apparatus, product, or process disclosed, or represents that its use would not infringe privately owned rights.
Reference herein to any specific commercial product, process, or service by trade name, trademark, manufacturer, or otherwise does not necessarily constitute or imply its endorsement, recommendation, or favoring by the United States government or Lawrence Livermore National Security, LLC.
The views and opinions of authors expressed herein do not necessarily state or reflect those of the United States government or Lawrence Livermore National Security, LLC, and shall not be used for advertising or product endorsement purposes.

\bibliographystyle{abbrv}
\bibliography{bibliography}

\begin{thebibliography}{10}

\bibitem{MFEM}
{MFEM}: Modular finite element methods library.
\newblock \url{https://mfem.org}.

\bibitem{Alexander1977}
R.~Alexander.
\newblock Diagonally implicit {R}unge-{K}utta methods for stiff {O.D.E.'s}.
\newblock {\em {SIAM} Journal on Numerical Analysis}, 14(6):1006--1021, Dec.
  1977.

\bibitem{Anderson2019}
R.~Anderson, J.~Andrej, A.~Barker, J.~Bramwell, J.-S. Camier, J.~Cerveny,
  V.~Dobrev, Y.~Dudouit, A.~Fisher, T.~Kolev, W.~Pazner, M.~Stowell, V.~Tomov,
  J.~Dahm, D.~Medina, and S.~Zampini.
\newblock {MFEM}: a modular finite element library, 11 2019.
\newblock arXiv preprint 1911.09220 (submitted for publication).

\bibitem{Bassi1997}
F.~Bassi and S.~Rebay.
\newblock A high-order accurate discontinuous finite element method for the
  numerical solution of the compressible {N}avier-{S}tokes equations.
\newblock {\em Journal of Computational Physics}, 131(2):267--279, 1997.

\bibitem{Bastian2019}
P.~Bastian, E.~H. M{\"u}ller, S.~M{\"u}thing, and M.~Piatkowski.
\newblock Matrix-free multigrid block-preconditioners for higher order
  discontinuous galerkin discretisations.
\newblock {\em Journal of Computational Physics}, 394:417--439, Oct 2019.

\bibitem{Bell1989}
J.~B. Bell, P.~Colella, and H.~M. Glaz.
\newblock A second-order projection method for the incompressible
  {N}avier-{S}tokes equations.
\newblock {\em Journal of Computational Physics}, 85(2):257--283, Dec. 1989.

\bibitem{Bello-Maldonado2019}
P.~D. Bello-Maldonado and P.~F. Fischer.
\newblock Scalable low-order finite element preconditioners for high-order
  spectral element {P}oisson solvers.
\newblock {\em SIAM Journal on Scientific Computing}, 41(5):S2--S18, 2019.

\bibitem{Benzi2005}
M.~Benzi, G.~H. Golub, and J.~Liesen.
\newblock Numerical solution of saddle point problems.
\newblock {\em Acta Numerica}, 14:1--137, Apr. 2005.

\bibitem{Brachet1984}
M.~E. Brachet, D.~Meiron, S.~Orszag, B.~Nickel, R.~Morf, and U.~Frisch.
\newblock The {T}aylor-{G}reen vortex and fully developed turbulence.
\newblock {\em Journal of Statistical Physics}, 34(5-6):1049--1063, Mar. 1984.

\bibitem{Brachet1983}
M.~E. Brachet, D.~I. Meiron, S.~A. Orszag, B.~G. Nickel, R.~H. Morf, and
  U.~Frisch.
\newblock Small-scale structure of the {T}aylor-{G}reen vortex.
\newblock {\em Journal of Fluid Mechanics}, 130(-1):411, May 1983.

\bibitem{Brenan1995}
K.~E. Brenan, S.~L. Campbell, and L.~R. Petzold.
\newblock {\em Numerical solution of initial-value problems in
  differential-algebraic equations}.
\newblock Society for Industrial and Applied Mathematics, Jan. 1995.

\bibitem{Brezzi1991}
F.~Brezzi and R.~S. Falk.
\newblock Stability of higher-order {H}ood--{T}aylor methods.
\newblock {\em {SIAM} Journal on Numerical Analysis}, 28(3):581--590, 1991.

\bibitem{Brix2014}
K.~Brix, C.~Canuto, and W.~Dahmen.
\newblock Nested dyadic grids associated with {L}egendre-{G}auss-{L}obatto
  grids.
\newblock {\em Numerische Mathematik}, 131(2):205--239, Dec. 2014.

\bibitem{Brown2001}
D.~L. Brown, R.~Cortez, and M.~L. Minion.
\newblock Accurate projection methods for the incompressible {N}avier-{S}tokes
  equations.
\newblock {\em Journal of Computational Physics}, 168(2):464--499, Apr. 2001.

\bibitem{Brown2018}
J.~Brown, A.~Abdelfata, J.-S. Camier, V.~Dobrev, J.~Dongarra, P.~Fischer,
  A.~Fisher, Y.~Dudouit, A.~Haidar, K.~Kamran, T.~Kalev, M.~Min, T.~Ratnayaka,
  M.~Shephard, C.~Smith, S.~Tomov, V.~Tomov, and T.~Warburton.
\newblock {CEED} {ECP} milestone report: Public release of {CEED} 1.0.
\newblock Technical Report WBS 2.2.6.06, CEED-MS13, U.S.\ Department of Energy,
  Mar. 2018.

\bibitem{Canuto1994}
C.~Canuto.
\newblock Stabilization of spectral methods by finite element bubble functions.
\newblock {\em Computer Methods in Applied Mechanics and Engineering},
  116(1-4):13--26, Jan. 1994.

\bibitem{Canuto2010}
C.~Canuto, P.~Gervasio, and A.~Quarteroni.
\newblock Finite-element preconditioning of {G}-{NI} spectral methods.
\newblock {\em {SIAM} Journal on Scientific Computing}, 31(6):4422--4451, Jan.
  2010.

\bibitem{Canuto1985}
C.~Canuto and A.~Quarteroni.
\newblock Preconditioned minimal residual methods for {C}hebyshev spectral
  calculations.
\newblock {\em Journal of Computational Physics}, 60(2):315--337, Sept. 1985.

\bibitem{Canuto2007}
C.~Canuto, A.~Quarteroni, M.~Y. Hussaini, and T.~A. Zang.
\newblock {\em Spectral methods}.
\newblock Springer Berlin Heidelberg, 2007.

\bibitem{Chorin1967}
A.~J. Chorin.
\newblock A numerical method for solving incompressible viscous flow problems.
\newblock {\em Journal of Computational Physics}, 2(1):12--26, Aug. 1967.

\bibitem{Chorin1968}
A.~J. Chorin.
\newblock Numerical solution of the {N}avier-{S}tokes equations.
\newblock {\em Mathematics of Computation}, 22(104):745--745, 1968.

\bibitem{Cockburn1989}
B.~Cockburn and C.-W. Shu.
\newblock {TVB} {R}unge-{K}utta local projection discontinuous {G}alerkin
  finite element method for conservation laws. {II}. {G}eneral framework.
\newblock {\em Mathematics of Computation}, 52(186):411--435, 1989.

\bibitem{Deville2002}
M.~O. Deville, P.~F. Fischer, E.~Mund, et~al.
\newblock {\em High-order methods for incompressible fluid flow}, volume~9.
\newblock Cambridge University Press, 2002.

\bibitem{Deville1990}
M.~O. Deville and E.~H. Mund.
\newblock Finite-element preconditioning for pseudospectral solutions of
  elliptic problems.
\newblock {\em {SIAM} Journal on Scientific and Statistical Computing},
  11(2):311--342, Mar. 1990.

\bibitem{Dobrev2017}
V.~Dobrev, J.~Dongarra, J.~Brown, P.~Fischer, A.~Haidar, I.~Karlin, T.~Kolev,
  M.~Min, T.~Moon, T.~Ratnayaka, S.~Tomov, and V.~Tomov.
\newblock {CEED} {ECP} milestone report: Identify initial kernels, bake-off
  problems (benchmarks) and miniapps.
\newblock Technical Report WBS 1.2.5.3.04, CEED-MS6, U.S.\ Department of
  Energy, 2017.

\bibitem{Dubief2000}
Y.~Dubief and F.~Delcayre.
\newblock On coherent-vortex identification in turbulence.
\newblock {\em Journal of Turbulence}, 1:N11, jan 2000.

\bibitem{Ewing1990}
R.~E. Ewing, R.~D. Lazarov, P.~Lu, and P.~S. Vassilevski.
\newblock Preconditioning indefinite systems arising from mixed finite element
  discretization of second-order elliptic problems.
\newblock In {\em Preconditioned Conjugate Gradient Methods}, pages 28--43,
  Berlin, Heidelberg, 1990. Springer Berlin Heidelberg.

\bibitem{Fehn2018}
N.~Fehn, W.~A. Wall, and M.~Kronbichler.
\newblock A matrix-free high-order discontinuous {G}alerkin compressible
  {N}avier-{S}tokes solver: A performance comparison of compressible and
  incompressible formulations for turbulent incompressible flows.
\newblock {\em International Journal for Numerical Methods in Fluids},
  89(3):71--102, Oct 2018.

\bibitem{Fischer1997}
P.~F. Fischer.
\newblock An overlapping {S}chwarz method for spectral element solution of the
  incompressible {N}avier--{S}tokes equations.
\newblock {\em Journal of Computational Physics}, 133(1):84 -- 101, 1997.

\bibitem{Fischer2005}
P.~F. Fischer and J.~W. Lottes.
\newblock Hybrid {S}chwarz-multigrid methods for the spectral element method:
  {E}xtensions to {N}avier-{S}tokes.
\newblock In {\em Domain Decomposition Methods in Science and Engineering},
  pages 35--49. Springer, 2005.

\bibitem{Franciolini2017}
M.~Franciolini, A.~Crivellini, and A.~Nigro.
\newblock On the efficiency of a matrix-free linearly implicit time integration
  strategy for high-order discontinuous {G}alerkin solutions of incompressible
  turbulent flows.
\newblock {\em Computers \& Fluids}, 159:276--294, Dec 2017.

\bibitem{Guermond2006}
J.-L. Guermond, P.~D. Minev, and J.~Shen.
\newblock An overview of projection methods for incompressible flows.
\newblock {\em Computer Methods in Applied Mechanics and Engineering},
  195:6011--6045, 2006.

\bibitem{Hairer1996}
E.~Hairer and G.~Wanner.
\newblock {\em Solving ordinary differential equations {II}}.
\newblock Springer Berlin Heidelberg, 1996.

\bibitem{Hunt1988}
J.~C.~R. Hunt, A.~A. Wray, and P.~Moin.
\newblock Eddies, streams, and convergence zones in turbulent flows.
\newblock Technical Report CTR-S88, Center for Turbulence Research, 1988.

\bibitem{John2006}
V.~John, G.~Matthies, and J.~Rang.
\newblock A comparison of time-discretization/linearization approaches for the
  incompressible {N}avier-{S}tokes equations.
\newblock {\em Computer Methods in Applied Mechanics and Engineering},
  195(44-47):5995--6010, Sept. 2006.

\bibitem{Karniadakis1991}
G.~E. Karniadakis, M.~Israeli, and S.~A. Orszag.
\newblock High-order splitting methods for the incompressible {N}avier-{S}tokes
  equations.
\newblock {\em Journal of Computational Physics}, 97(2):414--443, Dec. 1991.

\bibitem{Klockner2009}
A.~Kl{\"o}ckner, T.~Warburton, J.~Bridge, and J.~Hesthaven.
\newblock Nodal discontinuous {G}alerkin methods on graphics processors.
\newblock {\em Journal of Computational Physics}, 228(21):7863--7882, Nov.
  2009.

\bibitem{Klockner2011}
A.~Kl{\"o}ckner, T.~Warburton, and J.~S. Hesthaven.
\newblock Viscous shock capturing in a time-explicit discontinuous {G}alerkin
  method.
\newblock {\em Mathematical Modelling of Natural Phenomena}, 6(3):57--83, 2011.

\bibitem{Kovasznay1948}
L.~I.~G. Kovasznay.
\newblock Laminar flow behind a two-dimensional grid.
\newblock {\em Mathematical Proceedings of the Cambridge Philosophical
  Society}, 44(1):58--62, 1948.

\bibitem{Kronbichler2012}
M.~Kronbichler and K.~Kormann.
\newblock A generic interface for parallel cell-based finite element operator
  application.
\newblock {\em Computers \& Fluids}, 63:135 -- 147, 2012.

\bibitem{Kronbichler2019}
M.~Kronbichler and K.~Ljungkvist.
\newblock Multigrid for matrix-free high-order finite element computations on
  graphics processors.
\newblock {\em {ACM} Transactions on Parallel Computing}, 6(1):1--32, May 2019.

\bibitem{Lottes2005}
J.~W. Lottes and P.~F. Fischer.
\newblock Hybrid multigrid/{S}chwarz algorithms for the spectral element
  method.
\newblock {\em Journal of Scientific Computing}, 24(1):45--78, 2005.

\bibitem{Orszag1980}
S.~A. Orszag.
\newblock Spectral methods for problems in complex geometries.
\newblock {\em Journal of Computational Physics}, 37(1):70 -- 92, 1980.

\bibitem{Orszag1986}
S.~A. Orszag, M.~Israeli, and M.~O. Deville.
\newblock Boundary conditions for incompressible flows.
\newblock {\em Journal of Scientific Computing}, 1(1):75--111, 1986.

\bibitem{Patera1984}
A.~T. Patera.
\newblock A spectral element method for fluid dynamics: Laminar flow in a
  channel expansion.
\newblock {\em Journal of Computational Physics}, 54(3):468 -- 488, 1984.

\bibitem{Pazner2019Efficient}
W.~Pazner.
\newblock Efficient low-order refined preconditioners for high-order
  matrix-free continuous and discontinuous {G}alerkin methods, Aug. 2019.
\newblock arXiv preprint 1908.07071 (submitted for publication).

\bibitem{Pazner2019High}
W.~Pazner, M.~Franco, and P.-O. Persson.
\newblock High-order wall-resolved large eddy simulation of transonic buffet on
  the {OAT15A} airfoil.
\newblock In {\em AIAA Scitech 2019 Forum}, 2019.

\bibitem{Pazner2018Approximate}
W.~Pazner and P.-O. Persson.
\newblock Approximate tensor-product preconditioners for very high order
  discontinuous {G}alerkin methods.
\newblock {\em Journal of Computational Physics}, 354:344--369, 2018.

\bibitem{Pazner2018b}
W.~Pazner and P.-O. Persson.
\newblock Interior penalty tensor-product preconditioners for high-order
  discontinuous {G}alerkin discretizations.
\newblock In {\em 2018 {AIAA} Aerospace Sciences Meeting}. American Institute
  of Aeronautics and Astronautics, Jan. 2018.

\bibitem{Prohl1997}
A.~Prohl.
\newblock {\em Projection and quasi-compressibility methods for solving the
  incompressible {N}avier-{S}tokes equations}.
\newblock Vieweg+Teubner Verlag, 1997.

\bibitem{Rang2007}
J.~Rang.
\newblock Design of {DIRK} schemes for solving the {N}avier-{S}tokes equations.
\newblock {\em Informatik-Berichte der Technischen Universit{\"a}t
  Braunschweig}, 2007-02, 2007.

\bibitem{Rudi2015}
J.~Rudi, A.~C.~I. Malossi, T.~Isaac, G.~Stadler, M.~Gurnis, P.~W.~J. Staar,
  Y.~Ineichen, C.~Bekas, A.~Curioni, and O.~Ghattas.
\newblock An extreme-scale implicit solver for complex {PDE}s: Highly
  heterogeneous flow in {E}arth's mantle.
\newblock In {\em Proceedings of the International Conference for High
  Performance Computing, Networking, Storage and Analysis}, SC '15, pages
  5:1--5:12, New York, NY, USA, 2015. ACM.

\bibitem{Saad1993}
Y.~Saad.
\newblock A flexible inner-outer preconditioned {GMRES} algorithm.
\newblock {\em {SIAM} Journal on Scientific Computing}, 14(2):461--469, Mar.
  1993.

\bibitem{Saad2003}
Y.~Saad.
\newblock {\em Iterative methods for sparse linear systems}, volume~82.
\newblock {SIAM}, 2003.

\bibitem{Shen2011}
J.~Shen, T.~Tang, and L.-L. Wang.
\newblock {\em Spectral Methods}.
\newblock Springer Berlin Heidelberg, 2011.

\bibitem{Shu2005}
C.-W. Shu, W.-S. Don, D.~Gottlieb, O.~Schilling, and L.~Jameson.
\newblock Numerical convergence study of nearly incompressible, inviscid
  {T}aylor-{G}reen vortex flow.
\newblock {\em Journal of Scientific Computing}, 24(1):1--27, July 2005.

\bibitem{Swirydowicz2019}
K.~{\'S}wirydowicz, N.~Chalmers, A.~Karakus, and T.~Warburton.
\newblock Acceleration of tensor-product operations for high-order finite
  element methods.
\newblock {\em The International Journal of High Performance Computing
  Applications}, 33(4):735--757, 2019.

\bibitem{Teukolsky2015}
S.~A. Teukolsky.
\newblock Short note on the mass matrix for {G}auss-{L}obatto grid points.
\newblock {\em Journal of Computational Physics}, 283:408--413, Feb. 2015.

\bibitem{Tomboulides1997}
A.~G. Tomboulides, J.~C.~Y. Lee, and S.~A. Orszag.
\newblock Numerical simulation of low {M}ach number reactive flows.
\newblock {\em Journal of Scientific Computing}, 12(2):139--167, Jun 1997.

\bibitem{vanRees2011}
W.~M. van Rees, A.~Leonard, D.~Pullin, and P.~Koumoutsakos.
\newblock A comparison of vortex and pseudo-spectral methods for the simulation
  of periodic vortical flows at high {R}eynolds numbers.
\newblock {\em Journal of Computational Physics}, 230(8):2794--2805, Apr. 2011.

\bibitem{Vermeire2017}
B.~C. Vermeire, F.~D. Witherden, and P.~E. Vincent.
\newblock On the utility of {GPU} accelerated high-order methods for unsteady
  flow simulations: A comparison with industry-standard tools.
\newblock {\em Journal of Computational Physics}, 334:497--521, 2017.

\bibitem{Wang2013}
Z.~J. Wang, K.~Fidkowski, R.~Abgrall, F.~Bassi, D.~Caraeni, A.~Cary,
  H.~Deconinck, R.~Hartmann, K.~Hillewaert, H.~T. Huynh, et~al.
\newblock High-order {CFD} methods: current status and perspective.
\newblock {\em International Journal for Numerical Methods in Fluids},
  72(8):811--845, 2013.

\bibitem{Wilhelm2001}
D.~Wilhelm and L.~Kleiser.
\newblock Domain decomposition method and fast diagonalization solver for
  spectral element simulations.
\newblock {\em Computational Fluid Dynamics 2000}, pages 429--434, 2001.

\end{thebibliography}

\end{document}